\documentclass[11pt,a4paper]{amsart}
\usepackage[T1]{fontenc}
\usepackage{lmodern}
\usepackage{amsmath,amssymb,mathtools}
\usepackage[left=2.7cm, right=2.7cm, top=2.6cm, bottom=2.6cm]{geometry}

\usepackage{microtype,booktabs,array,enumitem}
\usepackage{needspace}
\usepackage{xcolor}
\usepackage[colorlinks=true,linkcolor=blue!50!black,citecolor=blue!50!black,urlcolor=blue!50!black]{hyperref}
\hypersetup{pdftitle={Multiplicity and non-hyperbolicity of closed geodesics on pinched Finsler spheres},pdfauthor={Guanjun Dong and Huagui Duan}}
\setlist[enumerate]{label=\textup{(\roman*)},leftmargin=2em}
\newtheorem{theorem}{Theorem}[section]
\newtheorem{lemma}[theorem]{Lemma}
\newtheorem{proposition}[theorem]{Proposition}
\newtheorem{corollary}[theorem]{Corollary}
\theoremstyle{definition}
\newtheorem{definition}[theorem]{Definition}
\theoremstyle{remark}
\newtheorem{remark}[theorem]{Remark}
\numberwithin{equation}{section}
\newcommand{\Q}{\mathbb Q}
\newcommand{\Z}{\mathbb Z}
\newcommand{\N}{\mathbb N}
\newcommand{\R}{\mathbb R}
\newcommand{\C}{\mathbb C}

\newcommand{\HH}{H^{S^1}}
\newcommand{\LL}{\Lambda S^n}
\newcommand{\LLz}{\Lambda^0 S^n}
\newcommand{\CC}{\overline C}
\newcommand{\ii}{\mathrm i}

\newcommand{\im}{\operatorname{im}}

\newcommand{\Vinf}{\mathcal V_\infty}
\newcommand{\Pcal}{\mathcal P}
\newcommand{\Nnh}{N_{\mathrm{nh}}}

\newcommand{\floor}[1]{\left\lfloor #1\right\rfloor}
\newcommand{\ceil}[1]{\left\lceil #1\right\rceil}

\title[Closed geodesics on pinched Finsler spheres]{Multiplicity and non-hyperbolicity of closed geodesics
on pinched Finsler spheres}
\author{Guanjun Dong}
\address{Guanjun Dong, School of Mathematical Sciences, Nankai University, Tianjin 300071, P. R. China}
\email{1120240031@mail.nankai.edu.cn}
\author{Huagui Duan}
\address{Huagui Duan, School of Mathematics, Jilin University, Changchun, Jinlin, 130012, P. R. China}
\email{duanhg@nankai.edu.cn}
\date{}
\begin{document}
\begin{abstract}
Let \(F\) be a Finsler metric on \(S^n\) with reversibility \(\lambda\) and
flag curvature \(K\) satisfying \((\lambda/(1+\lambda))^2<K\leq1\).
We prove that there exist at least \(n/2+1\) distinct prime closed geodesics when \(n\geq4\) is even. For \(S^4\), there are at least four distinct prime closed geodesics if it has at least two irrationally elliptic prime closed geodesics, and
if there are finitely many prime closed geodesics, there exist at least three non-hyperbolic prime
closed geodesics. Under the finiteness assumption, we further prove
that an odd-dimensional sphere \(S^{2k+1}\) has at least
\(k+1\) non-hyperbolic prime closed geodesics for \(k\geq2\).
\end{abstract}
\maketitle

\section{Introduction and main results}

\subsection{Introduction}
The multiplicity and stability of closed geodesics on Finsler spheres
are closely related to the Katok's examples and the Anosov conjecture \cite{Anosov75}.
Katok constructed irreversible Finsler metrics on spheres with
constant flag curvature and only finitely many prime closed geodesics
\cite{Katok73, Ziller83}.
These examples motivate the following conjecture:
\[
 N(S^n,F)\geq 2\floor{\frac{n+1}{2}},
\]
where $N(S^n,F)$ denotes the number of distinct prime closed geodesics
on \((S^n,F)\). Here a prime closed geodesic is one that is not a multiple
covering of another closed geodesic. Bangert and Long \cite{BangertLong10} proved the
existence of two distinct prime closed geodesics for every Finsler
$S^2$.

For a compact Finsler manifold, the reversibility is defined by
\[
 \lambda=\max\{F(-v):F(v)=1\}\geq1.
\]
Throughout this paper, the curvature assumption called the pinching condition is
\begin{equation}\label{eq:pinching}
 \left(\frac{\lambda}{1+\lambda}\right)^2<K\leq1,
\end{equation}
where \(K\) denotes the flag curvature.
The importance of this condition is that it gives uniform lower
bounds for the length, the Morse index and the mean index of every
closed geodesic. These comparison estimates were developed by
Rademacher \cite{Rad04,Rad07}. 

Without a non-degeneracy assumption, the conjecture remains difficult.
Under \eqref{eq:pinching}, Wang established multiplicity and stability
results in \cite{W08,W13,W19}. In particular, \cite[Theorem 1.2]{W19}
gives at least \(\floor{(n+1)/2}\) prime closed geodesics.
If the total number of prime closed geodesics is finite,
the results of \cite[Theorems 1.3 and 1.4]{W19} give at least
\(\floor{n/2}\) non-hyperbolic prime closed geodesics and at least
\(\floor{(n-1)/2}\) prime closed geodesics whose linearized
Poincar\'e maps have an eigenvalue \(\exp(\pi\ii\theta)\) with
\(\theta\notin\Q\), respectively.
Duan proved that there are at least three prime closed geodesics
for \(n\geq3\). If there are finitely many prime closed geodesics,
two of them can be chosen to be elliptic with an irrational unit-circle
eigenvalue \cite{Duan15,DuanEll15}. Moreover, \cite[Theorem 1.1]{Duan15}
gives three non-hyperbolic prime closed geodesics when \(n\geq6\).
For bumpy metrics, substantially stronger conclusions are available.
In particular, Wang proved the conjectured multiplicity under
\eqref{eq:pinching} \cite{W12}. Duan, Long and Wang proved that there exist at least n non-hyperbolic prime closed geodesics of positive index on an even-dimensional sphere, for bumpy irreversible metrics with finitely many prime closed geodesics \cite[Theorem~1.1]{DLW16}. For odd \(n\geq3\), their
Theorems~1.5 and~1.6 give at least \(n+1\) prime closed geodesics with even indices,
of which at least \(n-1\) are non-hyperbolic, under \eqref{eq:pinching}
and for bumpy metrics with finitely many prime closed geodesics.

A stronger pinching condition also leads to stronger conclusions without
a bumpy assumption. Recently, Duan and Xie proved the
existence of at least \(n\) prime closed geodesics for \(n\geq4\)
under the condition
\begin{equation}\label{eq:dx-pinching}
 \left(\frac{2n-3}{n-1}\right)^2
 \left(\frac{\lambda}{1+\lambda}\right)^2<K\leq1,
 \qquad \lambda<\frac{n-1}{n-2},
\end{equation}
and obtained the corresponding results on non-hyperbolicity and irrational eigenvalues
 \cite{DX26}. Their bound on the number of non-hyperbolic prime closed geodesics on \(S^4\)
under \eqref{eq:dx-pinching} is \(3\). The results below hold under \eqref{eq:pinching}, with no additional
restriction on \(\lambda\).

The equivariant method used in this paper has a long history in
closed geodesic theory. Rademacher developed a relative version of the Fadell--Rabinowitz index over \(H^*(BS^1)\) and the corresponding critical point theory in \cite[Sections 3, 5 and 6]{Rad94}. In equivariant
loop homology, Wang constructed the corresponding homology classes and proved a
strict ordering of their critical values, assuming that there are finitely many prime closed geodesics \cite[Propositions 2.1--2.3]{WExist12}. He also used
these classes to study non-hyperbolicity under additional metric
comparison and curvature hypotheses \cite[Theorem 1.4]{WExist12}.
Ginzburg and G\"urel established the
Lusternik--Schnirelmann theory and the construction of the so-called carriers on the
cotangent bundle via symplectic homology
\cite[Theorem 2.1 and Proposition 3.12]{GG20}.
Multiplicity and non-hyperbolicity on prequantization bundles are
studied in \cite{GGM18}. We give a direct
construction in equivariant loop homology in Section~\ref{sec:chains}. Within this framework, we can compare the contributions of the critical modules of geodesics with different iterations more intuitively.
One can also expect more applications of techniques from symplectic geometry to Finsler manifolds.

\subsection{Main results}
Denote \(\Pcal(S^n,F)\) by the set of distinct prime closed geodesics
and \(N(S^n,F)=\#\Pcal(S^n,F)\).
Our first result concerns even-dimensional spheres.

\begin{theorem}\label{thm:pinching}
Let \(n=2k\geq4\). Every Finsler metric \(F\) on \(S^n\) satisfying
\eqref{eq:pinching} corresponds to at least \(k+1\) distinct prime closed
geodesics.
\end{theorem}

\begin{remark}
Theorem~\ref{thm:pinching} raises Wang's general lower bound by one
in even dimensions in \cite{W19}. Its multiplicity conclusion on \(S^4\) is
already contained in \cite[Theorem~1.1]{Duan15}.
If the number of prime closed geodesics is finite, Wang's
Theorems~1.3 and~1.4 in \cite{W19} also guarantee at least \(k\)
non-hyperbolic prime closed geodesics.
\end{remark}

For a prime closed geodesic \(c\), let
\(P_c\in\operatorname{Sp}(2n-2)\) be its linearized Poincar\'e map.
The geodesic is hyperbolic if \(P_c\) has no eigenvalues on the unit
circle, and elliptic if all its eigenvalues lie on the unit circle.
Following \cite{Long02}, irrational ellipticity means that
\(\Omega^0(P_c)\) contains a \(\diamond\)-product of \(n-1\)
rotation matrices, whose angles are irrational multiples of \(\pi\).

\begin{theorem}\label{thm:four-conditional}
Let \((S^4,F)\) satisfy \eqref{eq:pinching}. If \(N(S^4,F)<\infty\), then there exist at least three non-hyperbolic prime closed geodesics. If there exist at least two irrationally elliptic prime closed geodesics, then \(N(S^4,F)\geq4\).
\end{theorem}

\begin{remark}
\cite[Lemma~3.1]{Duan15} gives the existence of two elliptic prime closed geodesics when \(N(S^4,F)<\infty\).
Here we need that the two elliptic prime closed geodesics are irrational.
\end{remark}

Let \(\Nnh(S^n,F)\) denote the number of non-hyperbolic prime closed
geodesics.

\begin{theorem}
\label{thm:odd-stability}
Let \(n=2k+1\geq5\), and let \((S^n,F)\) satisfy
\eqref{eq:pinching}. If \(N(S^n,F)<\infty\), then
\[
 \Nnh(S^n,F)\geq k+1=\frac{n+1}{2}.
\]
\end{theorem}

\begin{remark}
This raises the bound \(\floor{n/2}\) in
\cite[Theorem~1.3]{W19} by one in odd dimensions \(n\geq5\).
Theorem~\ref{thm:four-conditional} and \ref{thm:odd-stability} give three non-hyperbolic prime
closed geodesics for \(n=4,5\), which extend the conclusion of
\cite[Theorem~1.1]{Duan15} of three non-hyperbolic geodesics to these dimensions.
\end{remark}

\subsection{Outline of the proofs}
The difficulty in the degenerate case is that a single isolated
critical orbit can have nonzero local homology in several degrees.
Counting degrees alone therefore does not distinguish distinct iterations.
We use the carriers from Lemma~\ref{lem:chain}, together with the strict ordering of their critical values,
to obtain carriers with strictly increasing energies.
For \(J=1\), this fact is already
proved in \cite[Propositions 2.1--2.3]{WExist12}. For every \(J\), a formulation of these chains on the cotangent bundle is given in \cite[Proposition 3.12]{GG20}. We include a corresponding construction in
equivariant loop homology in Section~\ref{sec:chains}.

For \(n=2k\), suppose there are at most \(k\) prime closed geodesics.
The first \(k\) degrees of the chain for \(J=1\) lie below the Fadell--Rabinowitz indices of every
nonprime iteration. Their carriers therefore exhaust the set of prime closed
geodesics, and we can order them by length. The shortest prime closed geodesic
has nonzero critical module only in degree zero. A counting argument for visibility shows
that every prime closed geodesic is infinitely variationally visible, so we can identify
the ordering by length with the ordering by mean index.
Then, by the common index jump theorem, we choose a tuple with the required
divisibility, arbitrarily small errors of mean index and
\(d\mid N\). A truncated Euler characteristic identity then forces one central iteration to
contribute only in degree \(2N+d\), for the same tuple.
The degrees strictly between \(2N\) and \(2N+d\) are
carried by all the other central iterations. The degrees below
\(2N\) must then be carried by all the preceding iterations.
The resonance relation orders the preceding iterations in the reverse
order of the lengths of the prime closed geodesics. It forces the iterate of the shortest
prime closed geodesic to contribute in positive degree, which contradicts the invariance
of the critical module under the equal-nullity condition.

For Theorem~\ref{thm:four-conditional}, suppose that there are exactly
three prime closed geodesics on \(S^4\), and two of which are the irrationally
elliptic geodesics  (cf. \cite{Duan15}).
Then the common-index tuples force the third geodesic to have a rational rotation number
and a normal form consisting of two blocks at the eigenvalue one. The Morse inequalities
and the carriers of Fadell--Rabinowitz indices
then identify two infinitely variationally visible prime closed
geodesics. Their order by length in
the chain for \(J=1\) of Corollary~\ref{cor:ordered} is incompatible with the order forced by a large
common-index window .

In odd dimensions, the hypothetical minimal number of
non-hyperbolic prime closed geodesics forces the carrier of Fadell--Rabinowitz index \(2N\)
to be hyperbolic. By arguments similar to the above, its order by length and the sign of the
index-induced error then yield a contradiction.

\smallskip\noindent
\emph{Organization.}
Section~\ref{sec:prelim} recalls critical modules, Morse theory
and the relative Fadell--Rabinowitz index. Section~\ref{sec:index-theory} states
the index iteration theory for symplectic paths.
Section~\ref{sec:chains} constructs the ordered carriers and
establishes the estimates about the common index. Section~\ref{sec:even} proves
the three main results.

\medskip

\textbf{Acknowledgements.}  G. Dong was partially supported by NNSFC (12271268). H. Duan was partially supported by NNSFC (12631007, 12271268 and 12361141812) and Natural Science Foundation of Tianjin (No. 25JCZDJC01030).

\medskip

\textbf{AI disclosure.} ChatGPT was used for proofreading and language polishing, as well as for suggesting  the use of Lemma~\ref{lem:duan-top-propagation} to simplify the proof.

\section{Critical point theory and the Fadell--Rabinowitz index}
\label{sec:prelim}

We recall the variational results in the form used for closed
geodesics in \cite{BangertLong10,DX26,W08,W19}.
General references for the Morse-theoretic framework are
\cite{Chang93,MawhinWillem89}, and the setting of closed geodesics is
developed in \cite{Klingenberg78,Mercuri77,Rad92}.
Throughout, \(\N=\{1,2,\ldots\}\), \(\N_0=\N\cup\{0\}\), and
all homology and cohomology groups have rational coefficients.

\subsection{Local critical modules and the cyclic action}
An auxiliary Riemannian metric gives
\(\Lambda=H^1(S^1,S^n)\) its Hilbert manifold structure. The action
\[
 (s\cdot\gamma)(t)=\gamma(t+s),\qquad s,t\in S^1=\R/\Z,
\]
is continuous. The Finsler metric \(F\) defines an energy functional \(E\) and a length functional \(L\) on \(\Lambda\) by
\[
 E(\gamma)=\frac12\int_0^1F(\dot\gamma(t))^2\,dt,
 \qquad L(\gamma)=\int_0^1F(\dot\gamma(t))\,dt.
\]
Both functionals are invariant under the \(S^1\)-action. The functional \(E\) is \(C^{1,1}\) on \(\Lambda\) and satisfies the Palais--Smale condition (cf.~\cite{Mercuri77}). Thus we can apply the deformation theorems in \cite{{Chang93},{MawhinWillem89}}. The critical points of \(E\) of positive energies are precisely the closed geodesics \(c: S^1 \to M\) of the Finsler structure. If \(c \in \Lambda\) is a closed geodesic, then \(c\) is a regular curve, i.e., \(\dot{c}(t) \neq 0\) for all \(t \in S^1\), and this implies that the second differential \(E''(c)\) of \(E\) at \(c\) exists. As usual, we define the index \(i(c)\) of \(c\) as the maximal dimension of subspaces of \(T_c\Lambda\) on which \(E''(c)\) is negative definite, and the nullity \(\nu(c)\) of \(c\) so that \(\nu(c) + 1\) is the dimension of the null space of \(E''(c)\).

The iteration map is \(\phi_m(\gamma)(t)=\gamma(mt)\) and we write
\(c^m=\phi_m(c)\). If \(c\) is a prime closed geodesic, the stabilizer of \(c^m\)
is \(\Z_m\). One has
\begin{equation}\label{eq:energy-length}
 E(c)=\tfrac12L(c)^2,\qquad L(c^m)=mL(c),\qquad
 E(c^m)=m^2E(c).
\end{equation}
Let \(i(c^m)\) denote the Morse index of the
\(m\)-th iteration, and the average index is defined by
\(\hat i(c)=\lim_{m\to\infty}i(c^m)/m\).
Let \(\Lambda^0\) be the constant loops, and set
\[
 \Lambda^a=\{E\leq a\},\qquad \Lambda^{a-}=\{E<a\},
 \qquad \Lambda(c)=\{E<E(c)\}.
\]
We say that \(c\) satisfies \emph{(Iso)} if \(S^1\cdot c^m\) is an
isolated critical orbit for every \(m\in\N\).
Whenever the number of prime closed geodesics is finite, every prime closed geodesic
satisfies (Iso), and each bounded energy interval contains only
finitely many critical orbits.

For an isolated critical orbit, define its critical module by
\begin{equation}\label{eq:critical-module}
 \CC_q(E,c^m)=H_q\bigl((\Lambda(c^m)\cup S^1\cdot c^m)/S^1,
                          \Lambda(c^m)/S^1\bigr).
\end{equation}
Let \(U_{c^m}\) be a local negative disk and \(N_{c^m}\) a local
characteristic manifold in a slice transverse to the critical orbit.
Their dimensions are \(i(c^m)\) and \(\nu(c^m)\), respectively.
Use the notation
\[
 U_{c^m}^-=U_{c^m}\cap\Lambda(c^m),\qquad
 N_{c^m}^-=N_{c^m}\cap\Lambda(c^m),\qquad
 \mathcal H_s(c^m)=H_s(N_{c^m}^-\cup\{c^m\},N_{c^m}^-).
\]
For the action induced by a generator \(T\) of \(\Z_m\), let
\[
 H_*^{\pm\Z_m}=\{v:T_*v=\pm v\},\qquad
 K_s(c^m)=\dim\mathcal H_s(c^m),\qquad
 K_s^\pm(c^m)=\dim\mathcal H_s(c^m)^{\pm\Z_m}.
\]
Set \(\beta_m=(-1)^{i(c^m)-i(c)}\) and
\(k_s(c^m)=K_s^{\beta_m}(c^m)\).

\begin{proposition}\label{prop:splitting}
Let \(c\) be a prime closed geodesic satisfying (Iso). Then
\begin{equation}\label{eq:splitting-module}
 \CC_q(E,c^m)\cong
 \left(H_{i(c^m)}(U_{c^m}^-\cup\{c^m\},U_{c^m}^-)
       \otimes\mathcal H_{q-i(c^m)}(c^m)\right)^{+\Z_m}.
\end{equation}
Consequently
\[
 \CC_q(E,c^m)\cong\mathcal H_{q-i(c^m)}(c^m)^{\beta_m\Z_m},
 \qquad \dim\CC_q(E,c^m)=k_{q-i(c^m)}(c^m).
\]
In particular, when \(\nu(c^m)=0\),
\begin{equation}\label{eq:nondegenerate-module}
 \CC_q(E,c^m)=
 \begin{cases}
 \Q,&q=i(c^m)\text{ and }i(c^m)-i(c)\in2\Z,\\
 0,&\text{otherwise}.
 \end{cases}
\end{equation}
\end{proposition}

The underlying local construction of critical points goes back to
Gromoll and Meyer \cite{GM69a}. The displayed splitting formula
is in \cite[Satz 6.11]{Rad92}. The factor \(\beta_m\) records the orientation character of
the negative disk.

Clearly, the integers \(k_s(c^m)\) are equal to \(0\) when \(s < 0\) or \(s > \nu(c^m)\) and can take only values \(0\) or \(1\) when \(s = 0\) or \(s = \nu(c^m)\).

\begin{proposition}[{cf.~\cite[Satz 6.13]{Rad92}}]\label{prop:endpoints}
Let \(c\) be a prime cloed geodesic satisfying (Iso).
The groups \(\mathcal H_s(c^m)\) vanish outside \(0\leq s\leq \nu(c^m)\).
For any \(m\in \mathbb N\), we have:
\begin{enumerate}
\item If \(k_0(c^m)=1\), then \(k_s(c^m)=0\) for \(1\leq s\leq \nu(c^m)\).
\item If \(k_{\nu(c^m)}(c^m)=1\), then \(k_s(c^m)=0\) for \(0\leq s<\nu(c^m)\).
\item If \(k_s(c^m)>0\) for some \(0<s<\nu(c^m)\), then
\(k_0(c^m)=k_{\nu(c^m)}(c^m)=0\).
\item In particular, if \(\nu(c^m)\leq2\), then only one of the integers \(k_s(c^m)\) can be non-zero.
\end{enumerate}
\end{proposition}

\begin{proposition}
\label{prop:characteristic-iteration}
Suppose \(c\) is a prime closed geodesic and satisfies (Iso). For any positive integers \(a\) and \(b\) satisfying \(a\mid b\), we have
\begin{enumerate}
\item If \(\nu(c^a)=\nu(c^b)\), then
\(\mathcal H_s(c^a)\cong\mathcal H_s(c^b)\) for every \(s\).
\item If \(K_{\nu(c^b)}(c^b)=1\), then
\(K_{\nu(c^a)}(c^a)=1\), and
\[
K_s(c^b) = 0, \quad 0 \leq s \leq \nu(c^b) - 1.
\]
\end{enumerate}
\end{proposition}

\begin{proof}
Part (i) is the isomorphism in
\cite[Proposition 2.5(i)]{W08} by the formulation of the
critical point theory of Gromoll and Meyer
\cite{GM69b}. For (ii), this is
\cite[Proposition 2.7]{W08}.
\end{proof}

\begin{lemma}[{cf.~\cite[Propositions~2.1, 2.5(i) and~2.6]{W08}}]
\label{lem:local}
For isolated iterated closed geodesics, the following statements hold.
\begin{enumerate}
\item \(\CC_q(E,c^m)=0\) outside \([i(c^m),i(c^m)+\nu(c^m)]\).
\item If \(\CC_{i(c^m)+\nu(c^m)}(E,c^m)\ne0\), it is
one-dimensional and every other local module of this iteration is zero.
\item If \(c\) is a prime closed geodesic and \(\CC_{i(c)}(E,c)\ne0\), then
\(K_0(c)=1\) and \(K_s(c)=0\) for \(s>0\).
\item If \(\nu(c^m)=\nu(c)\), then \(K_s(c^m)=K_s(c)\).
Under the hypothesis of (iii), this implies
\(\CC_{i(c^m)+s}(E,c^m)=0\) for every \(s>0\).
\end{enumerate}
\end{lemma}

\begin{proposition}[{cf.~\cite[Lemma~5.2]{W08}}]\label{def:period}
Let \(c\) be a prime closed geodesic satisfying (Iso), then there exists a minimal \(T(c) \in \mathbb{N}\) such that
\[
\begin{aligned}
    \nu(c^{p+T(c)}) &= \nu(c^p),
    \quad i(c^{p+T(c)})-i(c^p)\in 2\mathbb{Z}
    && \text{for all } p\in\mathbb{N}, \\
    k_s(c^{p+T(c)}) &= k_s(c^p)
    && \text{for all } p\in\mathbb{N},\ s\in\mathbb{Z}.
\end{aligned}
\]
\end{proposition}

The Euler characteristic \(\chi(c^m)\) of \(c^m\) is defined by
\begin{equation}\label{eq:euler}
 \chi(c^m)=\sum_q(-1)^q\dim\CC_q(E,c^m)
          =\sum_{s=0}^{\nu(c^m)}(-1)^{i(c^m)+s}k_s(c^m),
\end{equation}
and the average Euler characteristic \(\hat\chi(c)\) of \(c\) is defined by
\begin{equation}\label{eq:average-euler}
 \hat\chi(c)=\frac1{T(c)}\sum_{m=1}^{T(c)}\chi(c^m).
\end{equation}
Thus \(\chi\) is periodic and \(\hat\chi(c)\in\Q\).

\begin{theorem}[{cf.~\cite[Theorem~7.9]{Rad92} or \cite[Theorem~5.5]{W08}}]\label{thm:mean-euler}
Suppose \((S^n,F)\) has finitely many prime closed geodesics
\(c_1,\ldots,c_p\), all with positive mean indices
\(\hat i(c_j)\).
Then
\begin{equation}\label{eq:resonance}
 \sum_{j=1}^p\frac{\hat\chi(c_j)}{\hat i(c_j)}=B(n,1),
 \qquad B(n,1)=\begin{cases}
 -\dfrac{n}{2n-2},&n\text{ even},\\[3pt]
 \dfrac{n+1}{2n-2},&n\text{ odd}.
 \end{cases}
\end{equation}
\end{theorem}

Set
$\overline{\Lambda}^{0}
=\overline{\Lambda}^{0}M
=\{\text{constant point curves in }M\}\cong M$.
Let $(X,Y)$ be a space pair such that the Betti numbers
$b_i=b_i(X,Y)=\dim H_i(X,Y;\mathbb{Q})$
are finite for all $i\in\mathbb{Z}$.
As usual, the \emph{Poincar\'e series} of $(X,Y)$ is defined by
the formal power series
$P(X,Y)=\sum_{i=0}^{\infty}b_i t^i$.
We need the following results on Betti numbers.

\begin{theorem}[{cf.~\cite{Hingston84},
    \cite[Theorem~2.4 and Remark~2.5]{Rad89}}]
We have the Poincar\'e series of
$(\overline{\Lambda}S^n,\overline{\Lambda}^{0}S^n)$
as the following.
\begin{enumerate}[label=\textup{(\roman*)}, leftmargin=*]
    \item
    When $n=2k+1$ is odd,
    \[
    \begin{aligned}
        P(\overline{\Lambda}S^n,\overline{\Lambda}^{0}S^n)(t)
        &=
        t^{n-1}
        \left(
            \frac{1}{1-t^2}
            +\frac{t^{n-1}}{1-t^{n-1}}
        \right) \\
        &=
        t^{2k}
        \left(
            \frac{1}{1-t^2}
            +\frac{t^{2k}}{1-t^{2k}}
        \right).
    \end{aligned}
    \]
    Thus, for $q\in\mathbb{Z}$ and $l\in\mathbb{N}_0$, we have
    \[
    \begin{aligned}
        b_q
        &=
        b_q(\overline{\Lambda}S^n,\overline{\Lambda}^{0}S^n)
        =
        \operatorname{rank}
        H_q(\overline{\Lambda}S^n,\overline{\Lambda}^{0}S^n) \\
        &=
        \begin{cases}
            2,
            & \text{if } q\in
              \{4k+2l \mid l\equiv 0 \pmod{k}\}, \\[2pt]
            1,
            & \text{if } q\in
              \{2k\}\cup
              \{2k+2l \mid l\not\equiv 0 \pmod{k}\}, \\[2pt]
            0,
            & \text{otherwise}.
        \end{cases}
    \end{aligned}
    \]

    \item
    When $n=2k$ is even,
    \[
    \begin{aligned}
        P(\overline{\Lambda}S^n,\overline{\Lambda}^{0}S^n)(t)
        &=
        t^{n-1}
        \left(
            \frac{1}{1-t^2}
            +\frac{t^{n(m+1)-2}}{1-t^{n(m+1)-2}}
        \right)
        \frac{1-t^{nm}}{1-t^n} \\
        &=
        t^{2k-1}
        \left(
            \frac{1}{1-t^2}
            +\frac{t^{4k-2}}{1-t^{4k-2}}
        \right),
    \end{aligned}
    \]
    where $m=1$ by \cite[Theorem~2.4]{Rad89}.
    Thus, for $q\in\mathbb{Z}$ and $l\in\mathbb{N}_0$, we have
    \[
    \begin{aligned}
        b_q
        &=
        b_q(\overline{\Lambda}S^n,\overline{\Lambda}^{0}S^n)
        =
        \operatorname{rank}
        H_q(\overline{\Lambda}S^n,\overline{\Lambda}^{0}S^n) \\
        &=
        \begin{cases}
            2,
            & \text{if } q\in
              \{6k-3+2l \mid l\equiv 0 \pmod{2k-1}\}, \\[2pt]
            1,
            & \text{if } q\in
              \{2k-1\}\cup
              \{2k-1+2l \mid l\not\equiv 0 \pmod{2k-1}\}, \\[2pt]
            0,
            & \text{otherwise}.
        \end{cases}
    \end{aligned}
    \]
\end{enumerate}
\end{theorem}

 In particular, let \(d=n-1=3\) and choose \(d\mid N\) for sufficiently large \(N\). Then
\begin{equation}\label{eq:lowdim-window-betti}
 b_{2N-d}=b_{2N+d}=2,\qquad
 b_{2N-d+2}=b_{2N+d-2}=1.
\end{equation}

We have the following version of the Morse inequality.

\begin{theorem}[{cf.~\cite[Theorem~6.1]{Rad92}}]
Suppose that there exist only finitely many prime closed geodesics
$\{c_j\}_{1\leq j\leq p}$ on $(M,F)$, and $0\leq a<b\leq\infty$.
Define, for each $q\in\mathbb{Z}$,
\[
\begin{aligned}
    M_q(\overline{\Lambda}^{b},\overline{\Lambda}^{a})
    &=
    \sum_{\substack{
        1\leq j\leq p,\ a<E(c_j^m)<b
    }}
    \operatorname{rank}\overline{C}_q(E,c_j^m), \\
    b_q(\overline{\Lambda}^{b},\overline{\Lambda}^{a})
    &=
    \operatorname{rank}
    H_q(\overline{\Lambda}^{b},\overline{\Lambda}^{a}).
\end{aligned}
\]
Then there holds
\[
\begin{aligned}
    &M_q(\overline{\Lambda}^{b},\overline{\Lambda}^{a})
    -M_{q-1}(\overline{\Lambda}^{b},\overline{\Lambda}^{a})
    +\cdots
    +(-1)^q M_0(\overline{\Lambda}^{b},\overline{\Lambda}^{a}) \\
    &\qquad\geq
    b_q(\overline{\Lambda}^{b},\overline{\Lambda}^{a})
    -b_{q-1}(\overline{\Lambda}^{b},\overline{\Lambda}^{a})
    +\cdots
    +(-1)^q b_0(\overline{\Lambda}^{b},\overline{\Lambda}^{a}), \\
    &M_q(\overline{\Lambda}^{b},\overline{\Lambda}^{a})
    \geq
    b_q(\overline{\Lambda}^{b},\overline{\Lambda}^{a}).
\end{aligned}
\]
\end{theorem}

\subsection{The relative Fadell--Rabinowitz index and critical values}
For an \(S^1\)-space \(X\) and a closed \(S^1\)-invariant subspace \(A \subset X\),
\[
 H^*_{S^1}(X,A)=H^*(ES^1\times_{S^1}X,ES^1\times_{S^1}A).
\]
This is the Borel construction in \cite{Borel60}.
Let \(u\) generate \(H^2(BS^1)\), and let
\(f:X_{S^1}\to BS^1\) be the classifying map.
The module action is \(u\cdot z=f^*u\smile z\) for \(z \in H^*_{S^1}(X,A)\).
We use the convention
\begin{equation}\label{eq:fr-order}
 \operatorname{ord}_u(z)
 =\inf\{r\in\N_0:u^r\cdot z=0\},
 \qquad \inf\varnothing=\infty.
\end{equation}
In particular, \(\operatorname{ord}_u(0)=0\).
This is the relative cohomological index used in
\cite{FR78,Rad94}. Specifying its value at zero fixes the indexing
of the critical levels below.

By \cite[Proposition 3.1]{Rad94}, there exists a class
\[
 z\in H^{n+1}_{S^1}(\Lambda,\Lambda^0),\qquad
 \operatorname{ord}_u(z)=\infty.
\]
For the inclusion \(j_a:(\Lambda^a,\Lambda^0)\hookrightarrow
(\Lambda,\Lambda^0)\), set
\begin{equation}\label{eq:fr-levels}
 d_z(a)=\operatorname{ord}_u(j_a^*z),\qquad
 \kappa_s=\inf\{a\geq0:d_z(a)\geq s\},\quad s\in\N.
\end{equation}
The function \(d_z\) is nondecreasing and right-continuous which we call it Fadell--Rabinowitz index (FR index).
Its discontinuities occur at critical values of \(E\), and a jump
of size at least two forces infinitely many prime closed geodesics,
cf. \cite[Section 5]{Rad94}.

\begin{theorem}
\label{thm:fr-carriers}
Suppose \((S^n,F)\) has only finitely many prime closed geodesics.
Then every \(\kappa_s\) is a positive critical value, and
\[
 0<\kappa_1<\kappa_2<\cdots,\qquad \kappa_s\longrightarrow\infty.
\]
For every \(s\in\N\), some prime closed geodesic \(c\) and iteration \(m\in\N\)
satisfy
\begin{equation}\label{eq:fr-carrier}
 E(c^m)=\kappa_s,\qquad
 \CC_{2s+n-1}(E,c^m)\ne0.
\end{equation}
\end{theorem}

The degree of the \(c^m\) is \(\deg z+2(s-1)=2s+n-1\).
The statements of critical value are in
\cite[Lemmas 2.3 and 2.4]{W13}.
The \(\kappa_s\) tends to infinity because a bounded energy interval
contains only finitely many critical orbits. We call such an iterate \(c^m\) a carrier of the critical value \(\kappa_s\) as in Section~\ref{sec:chains}.

\begin{definition}\label{def:visible}
A prime closed geodesic \(c\) is \((m,s)\)-variationally visible if
\eqref{eq:fr-carrier} holds. It is infinitely variationally visible
if this happens for infinitely many pairs \((m,s)\).
The set of such prime closed geodesics is denoted by \(\Vinf\).
\end{definition}

\begin{theorem}[{Resonance identity in \cite[Theorem~2.6]{W13}}]
\label{thm:visible-resonance}
When there are finitely many prime closed geodesics, the limit
\[
 \sigma=\lim_{s\to\infty}\frac{s}{\sqrt{2\kappa_s}}
\]
exists, and every \(c\in\Vinf\) satisfies
\begin{equation}\label{eq:wang-visible-ratio}
 \frac{\hat i(c)}{L(c)}=2\sigma.
\end{equation}
\end{theorem}

For each prime closed geodesic outside \(\Vinf\), the set of its visibility pairs
is finite by definition. The union of these sets over finitely many
prime closed geodesics is finite. Then we can have the following.

\begin{corollary}\label{cor:visibility-cutoff}
Suppose that there are finitely many prime closed geodesics.
There is an integer \(s_0\) such that every prime closed geodesic that is
\((m,s)\)-variationally visible for some \(s\geq s_0\) belongs to
\(\Vinf\). If all mean indices are positive, all such carriers
satisfy \(\hat i(c)/L(c)=\rho\) for the same \(\rho>0\).
\end{corollary}

\section{Index iteration theory for closed geodesics}
\label{sec:index-theory}

In this section, we recall briefly the index theory for symplectic paths developed by Long
and his coworkers, and all the details can be found in \cite{Long02}.
Let \(d=n-1\) and
\begin{gather*}
 J_d=\begin{pmatrix}0&-I_d\\ I_d&0\end{pmatrix},\qquad
 \operatorname{Sp}(2d)=\{M\in GL(2d,\mathbb R):M^TJ_dM=J_d\},\\
 \mathcal P_\tau(2d)=\{\gamma\in C([0,\tau],\operatorname{Sp}(2d)):
                         \gamma(0)=I_{2d}\}.
\end{gather*}
For \(\omega\in\mathbb U=\{z\in\C:|z|=1\}\), let
\(\operatorname{Sp}(2d)^0_\omega=\{M:\det(M-\omega I)=0\}\)
and let \(\operatorname{Sp}(2d)^*_\omega\) be its complement.
For any \(M \in \mathrm{Sp}(2d)_{\omega}^{0}\), we define a co-orientation of \(\mathrm{Sp}(2d)_{\omega}^{0}\) at \(M\) by the positive direction \(\left. \frac{d}{dt} M e^{t\epsilon J_d} \right|_{t=0}\) of the path \(M e^{t\epsilon J_d}\) with \(0 \le t \le 1\) and small \(\epsilon > 0\).

For any two continuous arcs \(\xi\) and \(\eta: [0, \tau] \to \mathrm{Sp}(2n)\) with \(\xi(\tau) = \eta(0)\), we define
\[
\eta * \xi(t) =
\begin{cases}
\xi(2t) & \text{if } 0 \le t \le \tau/2, \\
\eta(2t - \tau) & \text{if } \tau/2 \le t \le \tau.
\end{cases}
\]
A special path \(\xi_d\) is defined by
\[
\xi_d(t) =
\begin{pmatrix}
2 - \frac{t}{\tau} & 0 \\
0 & \left(2 - \frac{t}{\tau}\right)^{-1}
\end{pmatrix}^{\diamond d}
\quad \text{for } 0 \le t \le \tau.
\]
The \(\omega\)-index in \cite{Long02} for a path with
nondegenerate endpoint is the homotopy intersection number
\[
 i_\omega(\gamma)=[\operatorname{Sp}(2d)^0_\omega:
                         \gamma*\xi_d],\qquad
 \xi_d(t)=\operatorname{diag}\bigl((2-t/\tau)I_d,
                            (2-t/\tau)^{-1}I_d\bigr).
\]
For a degenerate endpoint, the
index is extended by
\[
 i_\omega(\gamma)=\sup_{\gamma \in U}
 \inf\{i_\omega(\eta):\eta\in U,
                       \eta(\tau)\in\operatorname{Sp}(2d)^*_\omega\},
 \qquad \nu_\omega(\gamma)=\dim_\C\ker_\C(\gamma(\tau)-\omega I).
\]

For \(m\in \mathbb N\), the iterated path is
\[
 \gamma^m(t)=\gamma(t-j\tau)\gamma(\tau)^j,
 \qquad j\tau\leq t\leq(j+1)\tau,\quad0\leq j<m.
\]
Set \(i(\gamma,m)=i_1(\gamma^m)\),
\(\nu(\gamma,m)=\nu_1(\gamma^m)\), and
\(\hat i(\gamma)=\lim_{m\to\infty}i(\gamma,m)/m\).
By \cite[Chapter 12]{Long02}, a
closed geodesic on an orientable Finsler manifold has a corresponding
symplectic path with endpoint \(P_c\) such that
\begin{equation}\label{eq:geodesic-index-convention}
 i(c^m)=i(\gamma_c,m),\qquad
 \nu(c^m)=\nu(\gamma_c,m)=\dim\ker(P_c^m-I),\qquad
 \hat i(c)=\hat i(\gamma_c).
\end{equation}
The relation between the geodesic Morse index and the symplectic
path index is discussed in \cite{Liu05,LiuLong02}, and index parity
is studied in \cite{Wilking01}.

For \(M=\gamma(\tau)\), the splitting numbers are defined by
\begin{equation}\label{eq:splitting-definition}
 S_M^\pm(\omega)=\lim_{\epsilon\to0^+}
          \bigl(i_{\omega\exp(\pm\ii\epsilon)}(\gamma)-i_\omega(\gamma)\bigr).
\end{equation}
They are nonnegative integers independent of the path used to
reach \(M\). They vanish off the spectrum and are additive under
symplectic direct sums. They are constant on the component
\(\Omega^0(M)\) containing \(M\) of the set
\[
 \Omega(M)=\{N:\sigma(N)\cap\mathbb U=\sigma(M)\cap\mathbb U,
 \ \nu_\omega(N)=\nu_\omega(M)\text{ for every }
                    \omega\in\sigma(M)\cap\mathbb U\}.
\]
In particular,\vspace{-2pt}
\begin{equation}\label{eq:splitting-basic}
 S^+_{N_1(1,b)}(1)=\begin{cases}1,&b\geq0,\\0,&b<0,\end{cases}
\end{equation}
which are shown in \cite{Long00} and \cite[Section 9.1]{Long02}.

The following symplectic matrices were introduced as basic normal forms.
For matrices \(M_j=\left(\begin{smallmatrix}A_j&B_j\\C_j&D_j\end{smallmatrix}\right)\),
their symplectic direct sum is defined by
\[
 M_1\diamond M_2=
 \begin{pmatrix}A_1&0&B_1&0\\0&A_2&0&B_2\\
                 C_1&0&D_1&0\\0&C_2&0&D_2\end{pmatrix}.
\]
The basic blocks used below are
\begin{gather*}
 D(\lambda)=\begin{pmatrix}\lambda&0\\0&\lambda^{-1}\end{pmatrix},
 \quad\lambda=\pm2,\qquad
 N_1(\lambda,b)=\begin{pmatrix}\lambda&b\\0&\lambda\end{pmatrix},
 \quad\lambda=\pm1,\ b\in\{-1,0,1\},\\
 R(\theta)=\begin{pmatrix}\cos\theta&-\sin\theta\\
                          \sin\theta&\cos\theta\end{pmatrix},
 \qquad N_2(e^{\ii\theta},B)=
       \begin{pmatrix}R(\theta)&B\\0&R(\theta)\end{pmatrix},
 \quad \theta\in(0,\pi)\cup(\pi,2\pi),
\end{gather*}
where \(B=(b_{ij})\) with \(b_{12}\ne b_{21}\) satisfying the
symplectic compatibility condition. The \(N_2(e^{\ii\theta},B)\) is called trivial
if \((b_{12}-b_{21})\sin\theta>0\), and nontrivial if this
quantity is negative. Let
\(\varphi(x)=\ceil{x}-\floor{x}\),
\(\{x\}=x-\floor{x}\).

The following theorem contains the precise index iteration formulae for symplectic paths, which is due to
Y. Long (cf. \cite[Theorem 8.3.1 and Corollary 8.3.2]{Long02}).

\begin{theorem}\label{thm:long-iteration}
Let \(\gamma\in\mathcal P_\tau(2d)\). Then there exists a path \(f \in C([0,1], \Omega^0(M))\) such that \(M:=f(0) = \gamma(\tau)\) and
\begin{align}
 f(1)={}&N_1(1,1)^{\diamond p_-}\diamond I_{2p_0}
       \diamond N_1(1,-1)^{\diamond p_+}\notag\\
 &\diamond N_1(-1,1)^{\diamond q_-}\diamond(-I_{2q_0})
       \diamond N_1(-1,-1)^{\diamond q_+}\notag\\
 &\diamond R(\theta_1)\diamond...\diamond R(\theta_r)
       \diamond N_2(e^{\ii\alpha_1},A_1) \diamond...\diamond N_2(e^{\ii\alpha_{r_*}},A_{r_*}) \\
 &\diamond N_2(e^{\ii\beta_1},B_1)\diamond...\diamond N_2(e^{\ii\beta_{r_0}},B_{r_0})
 \diamond M_0,\label{eq:full-normal-form}
\end{align}
where the \(N_2(e^{\ii\alpha_j},A_j)\) are nontrivial, the \(N_2(e^{\ii\beta_j},B_j)\)
are trivial, and \(\sigma(M_0)\cap\mathbb U=\varnothing\).
All multiplicities are nonnegative integers, and all angles belong
to \((0,\pi)\cup(\pi,2\pi)\). If \(M_0\in\operatorname{Sp}(2h)\),
then
\begin{equation}\label{eq:dimension}
 p_-+p_0+p_++q_-+q_0+q_++r+2r_*+2r_0+h=d.
\end{equation}
Set \(i_m=i(\gamma,m)\) and \(\nu_m=\nu(\gamma,m)\), one has
\begin{align}
 i_m={}&m(i_1+p_-+p_0-r)
       +2\sum_{\ell=1}^r\ceil{\frac{m\theta_\ell}{2\pi}}
       -r-p_--p_0\notag\\
 &-\frac{(1+(-1)^m)}{2}(q_0+q_+)+2\sum_{\ell=1}^{r_*}
             \left(\varphi\left(\frac{m\alpha_\ell}{2\pi}\right)-1\right),
             \label{eq:long-index}\\
 \nu_m={}&\nu_1+-\frac{(1+(-1)^m)}{2}(q_-+2q_0+q_+)
       +2\sum_{\ell=1}^{r}\left(1-\varphi\left(\frac{m\theta_\ell}{2\pi}\right)\right)
       \notag\\
 &+2\sum_{\ell=1}^{r_*}\left(1-\varphi\left(\frac{m\alpha_\ell}{2\pi}\right)\right)
       +2\sum_{\ell=1}^{r_0}\left(1-\varphi\left(\frac{m\beta_\ell}{2\pi}\right)\right),
       \label{eq:long-nullity}\\
 \hat i(\gamma)={}&i_1+p_-+p_0-r+\sum_{\ell=1}^r\frac{\theta_\ell}{\pi}.
       \label{eq:long-mean}
\end{align}
Moreover,
\begin{equation}\label{eq:splitting}
 S_M^+(1)=p_-+p_0,\qquad \nu_1=p_-+2p_0+p_+.
\end{equation}
We have \(i(\gamma,1)\) is odd if \(f(1)=N_1(1,1), I_2, N_1(-1,1), -I_2, N_1(-1,-1)\) and \(R(\theta)\); \(i(\gamma,1)\) is even if \(f(1)= N_1(1,-1)\) and \(N_2(e^{\ii\theta},B)\); \(i(\gamma,1)\) can be any integer if \(\sigma(f(1)) \cap \mathbb U = 0\).

\end{theorem}

\begin{theorem}\label{thm:successive}
Let \(M=\gamma(\tau)\), and let \(e(M)\) be the total algebraic
multiplicity of its unit-circle eigenvalues. Then for \(m\geq1\),
\begin{equation}\label{eq:successive-index}
 \nu_m-\frac{e(M)}2
 \leq i_{m+1}-i_m-i_1
 \leq\nu_1-\nu_{m+1}+\frac{e(M)}2.
\end{equation}
\end{theorem}

It was obtained in \cite[Theorem 2.2]{LongZhu02}.
In particular, if \(i_1\geq d\), then
\(i_{m+1}\geq i_m+\nu_m\), since \(e(M)\leq2d\).

%\begin{lemma}\label{lem:hyperbolic}
%If \(c\) is hyperbolic, then for every \(m\geq1\),
%\[
% \nu(c^m)=0,\qquad i(c^m)=m\,i(c),\qquad \hat i(c)=i(c)\in\Z.
%\]
%\end{lemma}
%\begin{proof}
%No eigenvalue of \(P_{c^m}\) is one, so every nullity is zero.
%By \(e(P_c)=0\), both bounds in
%\eqref{eq:successive-index} are zero. Induction gives the index
%identity, and division by \(m\) gives the mean index.
%\end{proof}

\begin{lemma}[{cf.~\cite[Theorem~10.1.2]{Long02}}]\label{lem:band}
For a closed geodesic on \(S^n\), with \(d=n-1\), one has
\begin{equation}\label{eq:band}
 m\hat i(c)-d\ \leq\ i(c^m)\ \leq\
 i(c^m)+\nu(c^m)\ \leq\ m\hat i(c)+d,\qquad r\geq1.
\end{equation}
\end{lemma}

Together the pinching assumption with \cite[Lemma~2]{Rad07}, we have the following.

\begin{lemma}
\label{lem:comparison}
Under \eqref{eq:pinching}, every prime closed geodesic \(c\) satisfies
\begin{equation}\label{eq:comparison}
 i(c^m)\geq md\quad(m\geq1),\qquad \hat i(c)>d.
\end{equation}
\end{lemma}

The common index jump theorem (cf. \cite[Theorem~4.3]{LongZhu02}) for symplectic paths has become one of the main tools in studying periodic orbits in Hamiltonian and symplectic dynamics. In 2016, the following enhanced common index jump theorem has been obtained by Duan, Long,and Wang in \cite{DLW16}.
In this paper, the iterations of
\(2m_j\) and \(2m_j-m\) are called by the central iterations and the preceding iterations respectively.

\begin{theorem}\label{thm:enhanced-cijt}
Let \(\gamma_j\in\mathcal P_{\tau_j}(2d)\), \(1\leq j\leq p\),
be finitely many paths with \(M_j=\gamma_j(\tau_j)\) and
\(\hat i(\gamma_j)>0\). Fix \(\bar m\in\N\).
There are arbitrarily large tuples \((N,m_1,\ldots,m_p)\) such
that, for every \(1\leq j\leq p\) and \(1\leq m\leq\bar m\),
\begin{align}
 \nu(\gamma_j,2m_j-m)&=\nu(\gamma_j,2m_j+m)=\nu(\gamma_j,m),
                       \label{eq:enhanced-nullity}\\
 i(\gamma_j,2m_j+m)&=2N+i(\gamma_j,m),\label{eq:enhanced-plus}\\
 i(\gamma_j,2m_j-m)&=2N-i(\gamma_j,m)
                         -2\bigl(S^+_{M_j}(1)+Q_j(m)\bigr),
                         \label{eq:enhanced-minus}\\
 i(\gamma_j,2m_j)&=2N-\bigl(S^+_{M_j}(1)+C(M_j)-2\Delta_j\bigr).
                         \label{eq:enhanced-center}
\end{align}
Here, with sufficiently small \(\delta>0\),
\begin{align*}
 C(M_j)&=\sum_{0<\theta<2\pi}S^-_{M_j}(e^{\ii\theta}),\qquad
 \Delta_j=\sum_{\substack{0<\theta<2\pi\\
                  0<\{m_j\theta/\pi\}<\delta}}S^-_{M_j}(e^{\ii\theta}),\\
 Q_j(m)&=\sum_{\substack{0<\theta<2\pi\\
             \{m_j\theta/\pi\}=0,\ \{m\theta/(2\pi)\}=0}}
                       S^-_{M_j}(e^{\ii\theta}).
\end{align*}
All sums are finite because splitting numbers vanish off the spectrum.

More precisely, choose an integer \(\bar M\) for which
\(\bar M\theta/\pi\in\Z\) whenever \(e^{\ii\theta}\in\sigma(M_j)\)
and \(\theta/\pi\in\Q\). The tuples can be chosen in the form
\begin{equation}\label{eq:enhanced-vertex}
 m_j=\bar M\left(\floor{\frac{N}{\bar M\hat i(\gamma_j)}}+\mathcal{X}_j\right),
 \qquad\mathcal{X}_j\in\{0,1\},\qquad
 \left|\left\{\frac{N}{\bar M\hat i(\gamma_j)}\right\}-\mathcal{X}_j\right|<\epsilon,
\end{equation}
for arbitrarily small \(\epsilon>0\).
One may require any fixed positive integer to divide \(N\),
and may enlarge \(\bar M\) to meet finitely many additional
divisibility requirements. These choices are simultaneous.
\end{theorem}

In particular, since \(0<\theta<2\pi\) implies \(\theta/(2\pi)\notin\Z\),
one has \(Q_j(1)=0\). Thus \eqref{eq:enhanced-minus} at \(m=1\)
contains exactly the correction \(2S^+_{M_j}(1)\).

\section{Ordered carriers and the common index of iterations}\label{sec:chains}

\subsection{The shift operator and ordered carriers}

We consider the relative equivariant homology
\[
 \HH_q(\LL,\LLz)
 =H_q(ES^1\times_{S^1}\LL,\ ES^1\times_{S^1}\LLz),\quad q\geq0.
\]
With the generator \(u\in H^2(BS^1)\) used in
Section~\ref{sec:prelim}, denote cap product with its pullback by
\(D=u\frown\), of degree \(-2\).
It is dual to the cup action used in the relative FR index:
\begin{equation}\label{eq:chain-pairing}
 \langle u\smile\alpha,w\rangle=\langle\alpha,Dw\rangle,
 \qquad
 \alpha\in H^{q-2}_{S^1}(\LL,\LLz),\quad
 w\in\HH_q(\LL,\LLz).
\end{equation}
Thus the notation describes the established equivariant module
structure \cite{Rad94,WExist12}. The construction below keeps
all classes and filtered maps in the setting of loop space.

For \(J=1\) in Corollary~\ref{cor:ordered} , it
gives the chain constructed by Wang in \cite[Proposition 2.1]{WExist12}.
The corresponding chains for every \(J\) appear in the symplectic
homology setting of cotangent bundle in
\cite[Proposition 3.12]{GG20}, where \(D\) is called the shift operator and the corresponding classes are called carriers. We show these constructions directly in the equivariant loop homology.

\begin{lemma}\label{lem:chain}
Let \(n\geq3\) and \(d=n-1\). For every \(J\geq1\), there are
classes satisfying
\begin{equation}\label{eq:chain}
 \begin{aligned}
 0\neq w_s^{(J)}&\in\HH_{(2J-1)d+2s}(\LL,\LLz),
                       &&0\leq s\leq d,\\
 Dw_{s+1}^{(J)}&=w_s^{(J)},&&0\leq s<d.
 \end{aligned}
\end{equation}
\end{lemma}

\begin{proof}
For \(J=1\), the construction in
\cite[Proposition 2.1]{WExist12} gives
\[
 w_d^{(1)}\in\HH_{3d}(\LL,\LLz),
 \qquad D^dw_d^{(1)}\neq0.
\]
Define
\[
 w_s^{(1)}=D^{d-s}w_d^{(1)},\qquad 0\leq s\leq d.
\]
Then
\[
 \deg w_s^{(1)}=3d-2(d-s)=d+2s,
 \qquad D^sw_s^{(1)}=D^dw_d^{(1)}\neq0.
\]
Thus \(w_s^{(1)}\neq0\) for every \(s\), and
\[
 Dw_{s+1}^{(1)}
 =D^{d-s}w_d^{(1)}=w_s^{(1)},\qquad 0\leq s<d.
\]

Let \(J\geq2\). By \cite[Proposition 3.1]{Rad94}, there is a class
\[
 z\in H^{d+2}_{S^1}(\LL,\LLz),
 \qquad u^az\neq0\quad\text{for every }a\geq0.
\]
Set
\[
 \ell=(J-1)d-1\geq0.
\]
Then we have the degree
\[
 \deg(u^{\ell+d}z)
 =d+2+2(\ell+d)
 =(2J+1)d.
\]
The universal coefficient theorem over \(\Q\) gives
\[
 H^{(2J+1)d}_{S^1}(\LL,\LLz)
 \cong
 \operatorname{Hom}_{\Q}
 \bigl(\HH_{(2J+1)d}(\LL,\LLz),\Q\bigr).
\]
Since \(u^{\ell+d}z\neq0\), choose and rescale
\(w_d^{(J)}\in\HH_{(2J+1)d}(\LL,\LLz)\) so that
\[
 \langle u^{\ell+d}z,w_d^{(J)}\rangle=1.
\]
For \(0\leq s\leq d\), define
\[
 w_s^{(J)}=D^{d-s}w_d^{(J)}.
\]
Then
\[
 \deg w_s^{(J)}
 =(2J+1)d-2(d-s)
 =(2J-1)d+2s.
\]
By \eqref{eq:chain-pairing}, we have
\[
 \begin{aligned}
 \langle u^{\ell+s}z,w_s^{(J)}\rangle
 &=\langle u^{\ell+s}z,D^{d-s}w_d^{(J)}\rangle\\
 &=\langle u^{\ell+d}z,w_d^{(J)}\rangle\\
 &=1.
 \end{aligned}
\]
Hence every \(w_s^{(J)}\) is nonzero. Finally,
\[
 Dw_{s+1}^{(J)}
 =D^{d-s}w_d^{(J)}
 =w_s^{(J)},\qquad 0\leq s<d.
\]
\end{proof}

We next recall the strict Lusternik--Schnirelmann inequality for the
shift operator \(D\). For the classes constructed at the chain for \(J=1\), strict ordering and the existence of local carriers are proved
in \cite[Propositions 2.2--2.3]{WExist12}. The general shift inequality
in equivariant Morse theory is \cite[Theorem 2.1 and Lemma 2.3]{GG20}, with its
relative version of filtered homology discussed in \cite[Remark 2.4]{GG20}.
We give the proof in relative equivariant homology of loop space.

\begin{lemma}
\label{lem:strict}
Suppose that there are finitely many prime closed geodesics. For a nonzero
homogeneous class \(w\in\HH_q(\LL,\LLz)\), define
\[
 c_F(w)=\inf\bigl\{a>0:w\in\im\bigl(
       \HH_q(\{E<a\},\LLz)\longrightarrow\HH_q(\LL,\LLz)
                              \bigr)\bigr\}.
\]
Then \(c_F(w)>0\) is a finite critical value, and there exist a critical orbit
\(S^1y\) at this value with \(\CC_q(E,y)\neq0\).
If \(Dw\neq0\), then
\begin{equation}\label{eq:strict}
 c_F(Dw)<c_F(w).
\end{equation}
\end{lemma}

\begin{proof}
Let \(\Lambda=\LLz\) and \(X_a=\{E<a\}\) for \(a>0\). Let
\[
 j_a:\HH_q(X_a,\Lambda)\longrightarrow\HH_q(\LL,\Lambda)
\]
be the map induced by the inclusion. On every relative group, \(D\) denotes the
cap product with the pullback of the same class
\(u\in H^2(BS^1)\). By naturality of the cap product,
\begin{equation}\label{eq:strict-naturality}
 j_a(Dz)=D(j_a z),
 \qquad z\in \HH_q(X_a,\Lambda).
\end{equation}
Let \(c_1,\ldots,c_p\) be the prime closed geodesics. Since
\(E(c_j^m)=m^2E(c_j)\), only finitely many of their iterates have
energy in a bounded interval. Consequently, every nonconstant
critical orbit is isolated, and there are only finitely many positive
critical values in each bounded interval.

\smallskip

{\bf Step 1}: {\it \(c_F(w)\) is a positive critical value.}

For sufficiently small \(a_0>0\), the deformation retraction of
small energy sublevels onto the constant loops gives
\[
H_*^{S^1}(X_a,\Lambda;\Q)=0,
\qquad 0<a\leq a_0,
\]
as used in the proof of
\cite[Proposition 2.2]{WExist12}.
Hence \(j_a=0\) for \(0<a\leq a_0\).
Since \(w\neq0\), we obtain
\[
c_F(w)\geq a_0>0.
\]
A relative cycle representing \(w\) in the Borel
construction consists of finitely many singular simplices, and its support is compact.
Since the energy \(E\) induces a continuous function on the
Borel construction, this support has bounded energy.
Thus \(w\in\im j_b\) for some finite \(b>0\), and
\[
0<c_F(w)<\infty.
\]
Set \(c=c_F(w)\).
The definition of \(c_F(w)\) and the inclusions
between sublevels imply
\[
w\in\im j_b
\qquad\text{for every }b>c.
\]
Suppose that \(c\) is not a critical value.
Then we can choose regular values
\(0<a<c<b\) such that \([a,b]\) contains no critical values.
The deformation lemma in \cite{Chang93,MawhinWillem89} gives an isomorphism
\[
H_q^{S^1}(X_a,\Lambda;\Q)
\longrightarrow
H_q^{S^1}(X_b,\Lambda;\Q).
\]
Consequently,
\[
w\in\im j_b=\im j_a,
\]
which contradicts \(a<c_F(w)\).
Therefore \(c_F(w)\) is a positive critical value.

\smallskip

{\bf Step 2:} {\it The Euler class vanishes near one critical orbit.}

Let \(\mathcal O=S^1y\) be a nonconstant critical orbit
with finite stabilizer \(\Gamma\).
Choose an \(S^1\)-invariant neighborhood \(U\) that
\(S^1\)-equivariantly retracts onto \(\mathcal O\),
as in the proof of
\cite[Proposition 2.2]{WExist12}.
Then
\[
U_{S^1}\simeq ES^1\times_{S^1}(S^1/\Gamma)
\cong B\Gamma.
\]
Since \(\Gamma\) is finite,
\begin{equation}\label{eq:local-u}
H^2(U_{S^1};\Q)=0,
\qquad u|_{U_{S^1}}=0.
\end{equation}
A local critical pair at \(\mathcal O\) lies in \(U\),
so the cap product with \(u\) vanishes on its
relative equivariant homology.

Choose regular values \(0<a<c<b\) such that
\(c=c_F(w)\) is the only critical value in \([a,b]\).
Let \(\mathcal O_j=S^1y_j\), \(1\leq j\leq s\),
be the critical orbits at energy \(c\) for some \(s>0\).
Equivariant deformation and excision give
\begin{equation}\label{eq:strict-localization}
\Phi_q:\HH_q(X_b,X_a)
\xrightarrow{\ \cong\ }
\bigoplus_{j=1}^{s}C_q^{S^1}(E,\mathcal O_j),
\end{equation}
as in
\cite[Proposition 2.3, equations (2.9)--(2.11)]{WExist12}.
By naturality of the cap product, the isomorphisms
in \eqref{eq:strict-localization} commute with \(D\).
Since \(D\) vanishes on each local critical module,
\[
 \Phi_{q-2}(Dz)=0
 \qquad\text{for every }z\in\HH_q(X_b,X_a).
\]
As \(\Phi_{q-2}\) is an isomorphism, we obtain
\begin{equation}\label{eq:strict-window-zero}
 D:\HH_q(X_b,X_a)\longrightarrow\HH_{q-2}(X_b,X_a),
 \qquad D=0.
\end{equation}
Since the stabilizers near \(\mathcal O_j\) are finite, the projection
from the Borel construction to the orbit space induces an isomorphism
on rational homology, as in the proof of \cite[Proposition~2.3, (2.11)]{WExist12}.
Thus
\[
C_q^{S^1}(E,\mathcal O_j)\cong\CC_q(E,y_j).
\]

\smallskip
{\bf Step 3:} \emph{There holds \(c_F(Dw)<c_F(w)\).}

The long exact sequence of the triple
\(\Lambda\subset X_a\subset X_b\) contains
\begin{equation}\label{eq:strict-triple}
\HH_q(X_a,\Lambda)\xrightarrow{i_q}\HH_q(X_b,\Lambda)
\xrightarrow{\pi_q}\HH_q(X_b,X_a).
\end{equation}
Since \(b>c\), choose
\(\widetilde w\in\HH_q(X_b,\Lambda)\) with
\(j_b\widetilde w=w\).
If \(\pi_q\widetilde w=0\), exactness gives
\(v\in\HH_q(X_a,\Lambda)\) such that
\(i_qv=\widetilde w\).
Then
\[
w=j_bi_qv=j_av,
\]
contrary to \(a<c_F(w)\).
Hence \(\pi_q\widetilde w\neq0\).
By \eqref{eq:strict-localization}, there is some \(j\)
such that
\[
\CC_q(E,y_j)\neq0,
\qquad E(y_j)=c_F(w).
\]

By naturality and \eqref{eq:strict-window-zero},
\[
\pi_{q-2}(D\widetilde w)
=D(\pi_q\widetilde w)=0.
\]
By applying \eqref{eq:strict-triple} in degree \(q-2\),
we obtain \(v'\in \HH_{q-2}(X_a,\Lambda)\) with
\[
i_{q-2}v'=D\widetilde w.
\]
Therefore
\[
j_av'
=j_bi_{q-2}v'
=j_b(D\widetilde w)
=D(j_b\widetilde w)
=Dw.
\]
If \(Dw\neq0\), it follows that
\[
c_F(Dw)\leq a<c_F(w).
\]
This proves \eqref{eq:strict}.
\end{proof}

Local vanishing of \(D\) allows one isolated orbit to have nonzero
local homology groups in two degrees differing by two. The strict inequality
concerns nonzero global classes related by \(D\). Such classes
have different critical values and therefore cannot have the same
critical orbit as a carrier.

\begin{corollary}\label{cor:ordered}
Suppose that \((S^n,F)\) corresponds to only finitely many prime closed geodesics.
Then for every \(J\geq1\), there are geodesics \(y_0,\ldots,y_d\) with
\begin{equation}\label{eq:ordered}
 \CC_{(2J-1)d+2s}(E,y_s)\neq0,\quad
 E(y_0)<E(y_1)<\cdots<E(y_d),\quad 0\leq s\leq d.
\end{equation}
In particular, no two of the geodesics \(y_0,\dots,y_d\) are iterates of the same prime closed geodesic.
\end{corollary}

\begin{remark}
For \(J=1\), the strict ordering and the argument of carriers are in
\cite[Propositions 2.2--2.3]{WExist12}. The proof of Proposition 2.3
there gives the nonzero local homology group as well as its index support.
The counterpart of this statement of cotangent bundle is
\cite[Corollary 3.14]{GG20}.
\end{remark}

\subsection{Variational visibility and resonance}

We apply Theorems~\ref{thm:fr-carriers} and
\ref{thm:visible-resonance} to obtain a common resonance relation for every prime closed geodesic
under the assumption that their number achieves the lower bound.

The following lemma can follow from~\cite[Theorem~6.4 and Corollary~6.5]{Rad94} and
Theorem~\ref{thm:visible-resonance}.
We include the argument needed for completeness.

\begin{lemma}\label{lem:visibility}
Suppose that there are finitely many prime closed geodesics and every mean index is larger than
\(d=n-1\). Then
\begin{equation}\label{eq:visibility-density}
 \frac12\leq\sum_{c_j\in\Vinf}\frac1{\hat i(c_j)},
 \qquad \#\Vinf>\frac d2.
\end{equation}
Thus, if the total number of prime closed geodesics is \(\ceil{n/2}\), all prime closed geodesics are
infinitely variationally visible. In that case, there is a common
\(\rho>0\) such that
\begin{equation}\label{eq:ratio}
 \hat i(c_j)=\rho L(c_j)\quad\text{for every prime closed geodesic }c_j.
\end{equation}
\end{lemma}

\begin{proof}
For each \(s\in\N\), choose a pair \((m(s),j(s))\) such that
\[
 E(c_{j(s)}^{m(s)})=\kappa_s,\qquad
 \CC_{2s+d}(E,c_{j(s)}^{m(s)})\ne0.
\]
A prime closed geodesic outside \(\Vinf\) has only finitely many variationally visible pairs.
Since there are only finitely many prime closed geodesics, there is a constant \(C_0\) such
that at most \(C_0\) of the chosen pairs involve geodesics outside \(\Vinf\).
For \(C>0\), set
\[
 N(C):=\#\{s\in\N:2s+d\leq C\}.
\]
Then \(N(C)=\frac C2+O(1)\).
If \(c_j^m\) corresponds to one of the degrees \(2s+d\), then
Lemma~\ref{lem:local}(i) and \eqref{eq:band} give
\[
 m\hat i(c_j)-d\leq i(c_j^m)\leq2s+d\leq C,
\]
and hence
\[
 m\leq\frac{C+d}{\hat i(c_j)}.
\]
Moreover, since the critical values \(\kappa_s\) are pairwise distinct,
a fixed iterate can correspond to at most one of them. Consequently,
\[
 N(C)-C_0=\frac C2+O(1)\leq\sum_{c_j\in\Vinf}\floor{\frac{C+d}{\hat i(c_j)}}.
\]
By dividing by \(C\) and letting \(C\to\infty\), we obtain
\[
 \frac12\leq\sum_{c_j\in\Vinf}\frac1{\hat i(c_j)}.
\]
Since \(\hat i(c_j)>d\) for every \(j\),
\[
 \frac12\leq\sum_{c_j\in\Vinf}\frac1{\hat i(c_j)}<\frac{\#\Vinf}{d}.
\]
Therefore
\[
 \#\Vinf>\frac d2,
\]
and hence
\[
 \#\Vinf\geq\ceil{\frac n2}.
\]
If the total number of prime closed geodesics is \(\ceil{n/2}\), every
prime closed geodesic belongs to \(\Vinf\).
Theorem~\ref{thm:visible-resonance} then yields
\[
 \frac{\hat i(c_j)}{L(c_j)}=2\sigma
\]
for every \(j\). Choosing \(\rho=2\sigma>0\) proves \eqref{eq:ratio}.
\end{proof}

By Theorem~\ref{thm:enhanced-cijt} with \(\bar m=1\), we can get the computations below directly.

\begin{lemma}\label{lem:tuple}
Let \(c_1,\dots,c_p\) be finitely many prime closed geodesics with \(\hat i(c_j)>0\), let \(T_j:=T(c_j)\) and \(\varepsilon>0\). There exist arbitrarily large integers \(N\) and
positive integers \(m_j\) such that
\begin{gather}
 d\mid N,\qquad T_j\mid 2m_j,\qquad
 \eta_j:=2m_j\hat i(c_j)-2N,\qquad |\eta_j|<\varepsilon,\label{eq:tuple-approx}\\
 \nu(c_j^{2m_j-1})=\nu(c_j),\label{eq:tuple-null}\\
 i(c_j^{2m_j-1})=2N-i(c_j)-2S^+_{P_{c_j}}(1),\label{eq:tuple-minus}\\
 i(c_j^{2m_j+1})=2N+i(c_j).\label{eq:tuple-plus}
\end{gather}
All these properties hold for the same tuple \((N,m_1,\ldots,m_p)\).
\end{lemma}

Following \cite[Theorem 4.2]{LongZhu02} and \cite{W19}, we can consider the sign conditions below.

\begin{lemma}\label{lem:signed-tuple}
Let \(c_1,\dots,c_p\) be as in Lemma~\ref{lem:tuple}. For every \(\varepsilon>0\), the tuples \((N,m_1,\dots,m_p)\) in Lemma~\ref{lem:tuple} can be chosen with \(N\) arbitrarily large so that
\[
\eta_j=0
\qquad\text{whenever }\hat i(c_j)\in\Q.
\]
For any fixed \(j_0\) such that \(\hat i(c_{j_0})\notin\Q\),
we can additionally require
\[
-\varepsilon<\eta_{j_0}<0.
\]
\end{lemma}

By the common index jump theorem, we can also get the supports of critical module in Lemma~\ref{lem:local} directly.

\begin{lemma}\label{lem:localization}
Suppose that \(i(c_j)\geq d\) and \(\hat i(c_j)>d\) for \(1\leq j\leq p\).
Choose
\[
0<\varepsilon<
\min\left\{1,\min_{1\leq j\leq p}(\hat i(c_j)-d)\right\}
\]
in Lemma~\ref{lem:tuple}.
Then for every \(1\leq j\leq p\), we have
\begin{equation}\label{eq:localization}
\left\{
\begin{aligned}
[i(c_j^m),\,i(c_j^m)+\nu(c_j^m)]
&\subset (-\infty,\,2N-d),
&&\text{for } 1\leq m\leq 2m_j-2,\\
[i(c_j^m),\,i(c_j^m)+\nu(c_j^m)]
&\subset (-\infty,\,2N),
&&\text{for } m=2m_j-1,\\
[i(c_j^m),\,i(c_j^m)+\nu(c_j^m)]
&\subset [2N-d,\,2N+d],
&&\text{for } m=2m_j,\\
[i(c_j^m),\,i(c_j^m)+\nu(c_j^m)]
&\subset [2N+d,\,+\infty),
&&\text{for } m=2m_j+1,\\
[i(c_j^m),\,i(c_j^m)+\nu(c_j^m)]
&\subset [2N+d,\,+\infty),
&&\text{for } m\geq 2m_j+2.
\end{aligned}
\right.
\end{equation}
\end{lemma}

\begin{proof}
The support is contained in \([i(c_j^m),\,i(c_j^m)+\nu(c_j^m)]\) by Lemma~\ref{lem:local}.
For \(m\leq2m_j-2\), its upper bound from \eqref{eq:band} is at most
\[
 2N+\eta_j-2\hat i(c_j)+d<2N-d.
\]
For \(m=2m_j-1\), it is at most
\(2N+\eta_j-\hat i(c_j)+d<2N\).
For \(m=2m_j\), Lemma~\ref{lem:band} and \(|\eta_j|<1\) give the desired interval.
For \(m=2m_j+1\), use the exact identity \eqref{eq:tuple-plus}
and \(i(c_j)\geq d\). For \(m\geq2m_j+2\), the index is at least
\(2N+\eta_j+2\hat i(c_j)-d>2N+d\).
\end{proof}

\begin{lemma}\label{lem:reverse}
Suppose that all prime closed geodesics are infinitely visible and are labeled by strictly
increasing lengths \(L(c_1)<\cdots<L(c_p)\).
Choose the same tuple with
\[
 2\varepsilon<\min_{i<j}(\hat i(c_j)-\hat i(c_i)).
\]
Then
\begin{equation}\label{eq:reverse}
 L(c_p^{2m_p-1})<\cdots<
 L(c_1^{2m_1-1}).
\end{equation}
\end{lemma}

\begin{proof}
By \eqref{eq:ratio}, \(\hat i(c_1)<\cdots<\hat i(c_p)\), and
\[
 \rho L(c_j^{2m_j-1})
   =(2m_j-1)\hat i(c_j)=2N-\hat i(c_j)+\eta_j.
\]
The difference between two error terms has absolute value below
\(2\varepsilon\), less than the mean index error. This gives
\eqref{eq:reverse}. Energy has the same order by
\eqref{eq:energy-length}.
\end{proof}

\begin{proposition}[{cf.~\cite[Lemma~4.1]{W19}}]\label{prop:top}
Suppose \(n\geq3\), there are finitely many prime closed geodesics,
\(i(c_j)\geq d\), and \(\hat i(c_j)>d\). Let
\((N,m_1,\ldots,m_p)\) be any sufficiently large tuple from
Lemma~\ref{lem:tuple} whose error \(\varepsilon\) satisfies
Lemma~\ref{lem:localization} and for which
\begin{equation}\label{eq:euler-error}
 \sum_j \frac{|\eta_j|\,|\hat\chi_j|}{\hat i(c_j)}<1.
\end{equation}
Then some \(j_*\) has
\begin{equation}\label{eq:upper-endpoint-support}
 \CC_{2N+d}(E,c_{j_*}^{2m_{j_*}})\neq0,\qquad
 \CC_q(E,c_{j_*}^{2m_{j_*}})=0\quad \text{for}~q\neq 2N+d.
\end{equation}
In particular, the conclusion holds for every
tuple satisfying either of the additional requirements in
Lemma~\ref{lem:signed-tuple}.
\end{proposition}

\section{Proof of main results}\label{sec:even}

\subsection{Multiplicity for even dimensions}

\begin{proof}[Proof of Theorem~\ref{thm:pinching}]
Let \(n=2k\geq4\) and \(d=2k-1\). Suppose for a contradiction
that there are \(p\leq k\) prime closed geodesics.

\smallskip

{\bf Step 1:} {\it Obtain a strict ordering of the mean indices.}

Apply Corollary~\ref{cor:ordered} with \(J=1\), and consider its first
\(k\) degrees:
\begin{equation}\label{eq:low-even}
 d,\ d+2,\ \ldots,\ d+2(k-1)=2d-1.
\end{equation}
Every degree in \eqref{eq:low-even} is below \(2d\).
By Lemma~\ref{lem:comparison}, every nonprime iteration
has index at least \(2d\), so such iterations from the degrees in \eqref{eq:low-even} can be excluded.
The corollary therefore provides \(k\) distinct prime closed geodesics with strictly increasing energies. Thus \(p=k\). Then we denote these prime closed geodesics by \(c_0,\ldots,c_{k-1}\) in the order of their
degrees. In particular, \(c_0\) is the shortest of all prime closed geodesics.

The nonzero critical module \(\CC_d(E,c_0)\), together with \(i(c_0)\geq d\),
implies \(i(c_0)=d\). Lemma~\ref{lem:local}(iii) then gives
\begin{equation}\label{eq:bottom}
 K_0(c_0)=1,\qquad K_s(c_0)=0\quad(s>0).
\end{equation}
Since \(p=k=\ceil{n/2}\), Lemma~\ref{lem:visibility} implies that
every prime closed geodesic is infinitely variationally visible.
Equation~\eqref{eq:ratio} therefore gives the same positive ratio
between mean index and length for all these prime closed geodesics. Their strict
ordering by length yields
\[
 \hat i(c_0)<\hat i(c_1)<\cdots<\hat i(c_{k-1}).
\]

\smallskip

{\bf Step 2:} {\it Determine the carriers in degrees strictly between \(2N\) and \(2N+d\).}

Choose an \(\varepsilon>0\) satisfying Lemmas~\ref{lem:reverse} and
\eqref{eq:euler-error}, then we can apply Lemma~\ref{lem:tuple} with \(\varepsilon\)
and sufficiently large \(N\).
For this same tuple, Proposition~\ref{prop:top} provides an index
\(j_*\) such that the only nonzero critical module of
\(c_{j_*}^{2m_{j_*}}\) occurs in degree \(2N+d\), as in
\eqref{eq:upper-endpoint-support}.
Since \(d\mid N\), Corollary~\ref{cor:ordered} holds for
\(J=N/d\). It gives \(2k\) carriers with strictly increasing
energies in the degrees
\begin{equation}\label{eq:high-even}
 2N-d,\ 2N-d+2,\ \ldots,\ 2N-1,\
 2N+1,\ \ldots,\ 2N+d.
\end{equation}

There are \(k-1\) degrees strictly between \(2N\) and \(2N+d\)
in \eqref{eq:high-even}:
\[
 2N+1,\ 2N+3,\ \ldots,\ 2N+d-2.
\]
By Lemma~\ref{lem:localization}, each corresponding carrier has the
form \(c_j^{2m_j}\). Their energies are strictly distinct, so the
\(k-1\) indices \(j\) are distinct. Also, no \(j\) equals \(j_*\) given by Proposition~\ref{prop:top} since
\(c_{j_*}^{2m_{j_*}}\) has zero critical module in each of these
degrees. Thus there are exactly \(k\) prime closed geodesics, and each iteration
\(c_j^{2m_j}\) with \(j\ne j_*\) corresponds to one of the displayed degrees.
It also follows that no iteration of the form \(c_j^{2m_j}\) can
correspond to a degree between \(2N-d\) and \(2N-1\).
For \(j=j_*\), its critical modules in that range vanish.
For \(j\ne j_*\), this iteration already corresponds to a carrier of higher degree. Since the carriers in \eqref{eq:ordered} have pairwise distinct critical values, it cannot also correspond to a carrier in this lower range.

\smallskip

{\bf Step 3:} {\it Determine the carriers in degrees less than \(2N\).}

There are \(k\) remaining degrees below \(2N\) in
\eqref{eq:high-even}:
\[
 2N-d,\ 2N-d+2,\ \ldots,\ 2N-1.
\]
Step 2 and Lemma~\ref{lem:localization} imply that every carrier
in these degrees belongs to the set
\(\{c_j^{2m_j-1}:0\leq j\leq k-1\}\).
The strict ordering of the critical values makes the \(k\) carriers
distinct. Since the set has exactly \(k\) elements, each element occurs once.
Lemma~\ref{lem:reverse} and Step 1 show
that \(c_0^{2m_0-1}\) has the largest energy in this set. It must
therefore correspond to the largest of the displayed degrees. Hence
\begin{equation}\label{eq:forced-nonzero}
 \CC_{2N-1}(E,c_0^{2m_0-1})\neq0.
\end{equation}

For the same tuple, Lemma~\ref{lem:tuple} gives
\[
 \nu(c_0^{2m_0-1})=\nu(c_0),\qquad
 i(c_0^{2m_0-1})=2N-d-2S^+_{P_{c_0}}(1).
\]
Equation~\eqref{eq:bottom} and Proposition~\ref{prop:characteristic-iteration} imply
\(K_s(c_0^{2m_0-1})=0\) for every \(s>0\).
By Proposition~\ref{prop:splitting}, we then have
\[
0\leq \CC_{2N-1}(E,c_0^{2m_0-1})=k_{s_0}(c_0^{2m_0-1}),
\]
where \(s_0\) equals
\begin{equation}\label{eq:offset-contradiction}
 (2N-1)-i(c_0^{2m_0-1})
 =d-1+2S^+_{P_{c_0}}(1)>0.
\end{equation}
Here \(d\geq3\) and the splitting number is nonnegative.
Since \(k_{s_0}(c_0^{2m_0-1})\leq K_{s_0}(c_0^{2m_0-1})=0\), this contradicts
\eqref{eq:forced-nonzero}. Thus \(p\leq k\) is impossible.
\end{proof}

\subsection{The four-dimensional result}

Recall that \(K_s(c^m)=\dim\mathcal H_s(c^m)\) and \(k_s(c^m)=\dim\mathcal H_s(c^m)^{\beta_m\Z_m}\) with
\(\beta_m=(-1)^{i(c^m)-i(c)}\). Thus
\begin{equation}\label{ineq:subspace}
 \dim\CC_{i(c^m)+s}(E,c^m)=k_s(c^m)\leq K_s(c^m).
\end{equation}

We begin with the two elliptic prime closed geodesics constructed by \cite{Duan15}. The following lemma follows from \cite[Lemmas~3.1 and 3.3]{Duan15} together with \cite[Propositions~2.5(i) and 2.7]{W08}. We include the argument for arbitrary iterations.

\begin{lemma}\label{lem:duan-top-propagation}
Suppose that \((S^n,F)\) satisfies \eqref{eq:pinching},
and that it has finitely many prime closed geodesics for \(n\geq3\).
There are two distinct elliptic prime closed geodesics \(c_1\) and \(c_2\), each possessing an eigenvalue
\(\exp(\pi\ii\theta)\) with \(\theta\notin\Q\), such that for every
\(m\geq1\),
\[
 K_s(c_j^m)=
 \begin{cases}
  1,&s=\nu(c_j^m),\\
  0,&s\ne\nu(c_j^m),
 \end{cases}
 \qquad j=1,2.
\]
Consequently,
\[
 \CC_q(E,c_j^m)=0\quad\text{if }q\ne i(c_j^m)+\nu(c_j^m),
 \qquad \sum_q\dim\CC_q(E,c_j^m)\leq1.
\]
\end{lemma}

\begin{proof}
By \cite[Lemmas 3.1 and 3.3]{Duan15}, there exist two
distinct elliptic prime closed geodesics \(c_1\) and \(c_2\)
with the irrational eigenvalues such that
\[
k_{\nu(c_j^{T_j})}(c_j^{T_j})=1,
\qquad T_j=T(c_j),\qquad j=1,2.
\]
Since
\[
1=k_{\nu(c_j^{T_j})}(c_j^{T_j})
\leq K_{\nu(c_j^{T_j})}(c_j^{T_j})
\leq1,
\]
we obtain
\[
K_{\nu(c_j^{T_j})}(c_j^{T_j})=1.
\]

Fix \(j\in\{1,2\}\), \(m\geq1\) and let \(L=\operatorname{lcm}(m,T_j)\).
Then the divisibility of \(T_j\mid L\) gives
\[
\ker(P_{c_j}^{T_j}-I)
\subseteq\ker(P_{c_j}^{L}-I).
\]
Since period \(T_j\) mentioned in Proposition~\ref{def:period} realizes the maximal nullity,
\[
\nu(c_j^{T_j})
\leq\nu(c_j^L)
\leq\max_{a\geq1}\nu(c_j^a)
=\nu(c_j^{T_j}).
\]
Thus \(\nu(c_j^L)=\nu(c_j^{T_j})\).
By \cite[Proposition 2.5(i)]{W08},
\[
K_{\nu(c_j^L)}(c_j^L)
=K_{\nu(c_j^{T_j})}(c_j^{T_j})
=1.
\]
If \(L>m\), then \(m\mid L\) and
\(\nu(c_j^m)\leq\nu(c_j^L)\).
Proposition 2.7 of \cite{W08} therefore implies
\[
K_{\nu(c_j^m)}(c_j^m)=1.
\]
The same equality already holds if \(L=m\).
The endpoint property in
Proposition~\ref{prop:endpoints} now gives
\[
K_s(c_j^m)=
\begin{cases}
1,&s=\nu(c_j^m),\\
0,&s\neq\nu(c_j^m).
\end{cases}
\]

Finally, Proposition~\ref{prop:splitting} gives
\[
\dim\CC_q(E,c_j^m)
=k_{q-i(c_j^m)}(c_j^m)
\leq K_{q-i(c_j^m)}(c_j^m).
\]
Consequently,
\[
\CC_q(E,c_j^m)=0
\quad\text{for }q\neq i(c_j^m)+\nu(c_j^m),
\qquad
\sum_q\dim\CC_q(E,c_j^m)\leq1.
\]
\end{proof}

\begin{proof}[Proof of Theorem~\ref{thm:four-conditional}]
Let \(d=n-1=3\). By Lemma~\ref{lem:duan-top-propagation}, there exist at least two elliptic
prime closed geodesics \(c_1\) and \(c_2\).

\smallskip

{\bf Step 1:} {\it Existence of at least three non-hyperbolic prime closed geodesics.}

Suppose that there are exactly two non-hyperbolic prime closed geodesics, and denote them by \(c_1\) and \(c_2\). Every other prime closed geodesic is then hyperbolic.
Lemma~\ref{lem:comparison} gives \(i(c)\geq d\) and \(\hat i(c)>d\)
for every prime closed geodesic \(c\). Choose a tuple from
Lemma~\ref{lem:tuple} with \(d\mid N\), large enough for
\eqref{eq:lowdim-window-betti} and compatible with
Lemma~\ref{lem:localization} and Proposition~\ref{prop:top}.
We can also choose the \(\varepsilon\) sufficiently small.

For a hyperbolic prime closed geodesic \(h\),
Theorem~\ref{thm:long-iteration} gives
\[
 \nu(h^m)=0,\qquad i(h^m)=mI,\qquad
 I=i(h)=\hat i(h)\in\Z,\quad I\geq d+1.
\]
The inequality \(|2m_hI-2N|<1\) and integrality imply
\(2m_hI=2N\) for some \(m_h\). If \(m\ne2m_h\), then
\[
 |i(h^m)-2N|=|m-2m_h|I\geq d+1.
\]
Thus the only possible contribution of \(h\) in
\([2N-d,2N+d]\) is in degree \(2N\).
In particular, there is no such contribution in the four degrees
\begin{equation}\label{eq:lowdim-four-degrees}
 2N-d,\quad 2N-d+2,\quad 2N+d-2,\quad 2N+d.
\end{equation}
By \eqref{eq:lowdim-window-betti}, the Betti numbers in these degrees are \(2,1,1,2\), respectively.

Using Proposition~\ref{prop:top} for the tuple above yields a geodesic with the
central iteration whose critical module is nozero in degree \(2N+d\) and
vanishes in any other degree. In fact, it belongs to one
of \(c_1\) and \(c_2\), and we can relabel these two geodesics so that it is \(c_1^{2m_1}\).
The Morse inequality in degree \(2N+d-2\) gives \(M_{2N+d-2}\geq1\).
By Lemma~\ref{lem:localization}, a closed geodesic
with iteration at most \(2m_j-1\) has all its critical modules
strictly below degree \(2N\), whereas one with iteration at least
\(2m_j+1\) has index at least \(2N+d\). Hence only a central
iteration can contribute in degree \(2N+d-2\). We can exclude
the hyperbolic ones and \(c_1^{2m_1}\) so
\begin{equation}\label{eq:second-elliptic-center}
 \CC_{2N+d-2}(E,c_2^{2m_2})\ne0.
\end{equation}
Lemma~\ref{lem:duan-top-propagation} now shows that
\(c_2^{2m_2}\) contributes in no other degree. In particular,
both of them with the central iterations have zero critical modules at
\(2N-d\) and \(2N-d+2\).

For the two lower degrees \(2N-d\) and \(2N-d+2\), we
exclude all hyperbolic iterates and the two elliptic geodesics with the central iterations. Lemma~\ref{lem:localization} shows that
the indices of iterations at most \(2m_j-2\) are strictly below
\(2N-d\), and those of iterations at least \(2m_j+1\) are no less than
\(2N+d\). The only remaining possibilities are therefore
\[
 c_1^{2m_1-1},\qquad c_2^{2m_2-1}.
\]
Each of these two iterates has critical module of total dimension at
most one by Lemma~\ref{lem:duan-top-propagation}. Consequently,
\[
 M_{2N-d}+M_{2N-d+2}\leq2.
\]
On the other hand, the Morse inequalities give
\[
 3=b_{2N-d}+b_{2N-d+2}
 \leq M_{2N-d}+M_{2N-d+2}\leq2,
\]
a contradiction.

\smallskip

{\bf Step 2:} {\it Existence of at least four prime closed geodesics. }

By Theorem~\ref{thm:pinching}, it suffices to rule out the case of exactly three
prime closed geodesics. Suppose that they are $c_1,c_2,c_3$, where
$c_1$ and $c_2$ are the irrationally elliptic geodesics in the hypothesis.
Lemma~\ref{lem:comparison} gives
\begin{equation}\label{eq:s4-basic-bounds}
 i(c_j^m)\geq3m,\qquad \hat i(c_j)>3
 \qquad(1\leq j\leq3,\ m\geq1).
\end{equation}
Each of $c_1$ and $c_2$ has odd index, and every iterate of them is nondegenerate
by Theorem~\ref{thm:long-iteration}. Hence
\[
 \CC_q(E,c_j^m)=
 \begin{cases}\Q,&q=i(c_j^m),\\0,&q\ne i(c_j^m),\end{cases}
 \qquad j=1,2.
\]

\smallskip
{\bf Step 2.1:} {\it Index estimates and local homology of \(c_3\).}

First, we consider the contribution of \(c_3\). We use the construction in \cite[Lemma~3.1]{Duan15} with the
fixed geodesics $c_1,c_2$, and denote the corresponding two tuples by
$(N,m_1,m_2,m_3)$ and $(N',m_1',m_2',m_3')$ which satisfy certain complementarity relations.
Set $\eta_j=2m_j\hat i(c_j)-2N$ and
$\eta_j'=2m_j'\hat i(c_j)-2N'$ for $j=1,2,3$.
The index formula in \cite[(3.12)--(3.14)]{Duan15} is
\[
 i(c_j^{2m_j})=2N-S^+_{P_{c_j}}(1)-C(P_{c_j})+2\Delta_j,
 \qquad
 \Delta_j=\sum_{0<\{m_j\theta/\pi\}<\delta}
                 S^-_{P_{c_j}}(e^{\ii\theta}),
\]
where $C(P_{c_j})=\sum_{0<\theta<2\pi}S^-_{P_{c_j}}(e^{\ii\theta})$.
For $j=1,2$, irrational ellipticity gives
$S^+_{P_{c_j}}(1)=0$ and $C(P_{c_j})=3$. The
equalities in \cite[(3.15)--(3.19)]{Duan15}
therefore imply
\[
 \Delta_1=\Delta_2'=3,\qquad
 \Delta_j+\Delta_j'=3\ (j=1,2),\qquad
 \Delta_2=\Delta_1'=0.
\]
Consequently,
\begin{equation}\label{eq:s4-pair-one}
\begin{aligned}
 i(c_1^{2m_1})&=2N+3,& i(c_2^{2m_2})&=2N-3,\\
 i(c_1^{2m_1'})&=2N'-3,& i(c_2^{2m_2'})&=2N'+3.
\end{aligned}
\end{equation}
For \(c_3\), \cite[(3.31)]{DuanEll15} gives the complementarity relation
\[
 \Delta_3+\Delta_3'
 =\sum_{\substack{0<\theta<2\pi\\
                  \theta/\pi\notin\Q}}
 S^-_{P_{c_3}}(e^{\ii\theta}).
\]

We impose the divisibility conditions of Lemma~\ref{lem:tuple} when considering these two tuples.
In Theorem~\ref{thm:enhanced-cijt}, we can choose the denominator
\(\bar M\) to be divisible by all \(T(c_j)\).
By \cite[Theorem 2.7 and (2.27)]{Duan15},
we also require \(3\mid N\).
So we can apply the construction to the two tuples, with arbitrarily small errors, as given by
\cite[Theorem 4.2]{LongZhu02}.
Thus, for each tuple,
\[
3\mid N,\qquad
T(c_j)\mid2m_j,\qquad
|\eta_j|=|2m_j\hat i(c_j)-2N|<\varepsilon,
\qquad j=1,2,3.
\]
The index identities and complementarity relations above
are preserved.
In what follows, each of the two tuples is chosen
with \(N\) sufficiently large and \(\varepsilon\)
sufficiently small for Lemma~\ref{lem:localization}.

Taking $J=N/3$ in Corollary~\ref{cor:ordered} gives a strict ordering
in the degrees
\begin{equation}\label{eq:s4-four-chain}
 2N-3,\qquad2N-1,\qquad2N+1,\qquad2N+3.
\end{equation}
Their Betti numbers are $2,1,1,2$.
Using the multiplicities of Theorem~\ref{thm:long-iteration}
and setting
\begin{equation}\label{eq:s4-G-def}
 G_3=i(c_3)+2S^+_{P_{c_3}}(1)-\nu(c_3)
     =i(c_3)+p_{3,-}-p_{3,+},
\end{equation}
we obtain
\[
 i(c_3^{2m_3-1})+\nu(c_3^{2m_3-1})=2N-G_3.
\]
If $p_{3,+}\leq2$, then $G_3\geq3-p_{3,+}\geq1$.
If $p_{3,+}=3$, all blocks of normal form are $N_1(1,-1)$,
so the prime index is even and $G_3=i(c_3)-3\geq1$.
Suppose $G_3\geq2$. For $j=1,2$,
\[
 i(c_j^{2m_j-1})=2N-i(c_j)\leq2N-3,
 \qquad \nu(c_j^{2m_j-1})=0,
\]
whereas the upper endpoint of the critical module for \(c_3\) is
$2N-G_3\leq2N-2$.
By \eqref{eq:s4-pair-one} and Lemma~\ref{lem:localization}, only
$c_3^{2m_3}$ can correspond to either degree $2N-1$ or $2N+1$ of
\eqref{eq:s4-four-chain}. Their critical values would both equal
$E(c_3^{2m_3})$, contradicting the strict ordering. Hence
\begin{equation}\label{eq:s4-G-one}
 G_3=1.
\end{equation}
Then the carrier in degree $2N+1$ is still $c_3^{2m_3}$. The carrier in degree
$2N-1$ must have smaller energy and is therefore
$c_3^{2m_3-1}$. Consequently,
\[
 \CC_{2N-1}(E,c_3^{2m_3-1})\ne0,\qquad
 i(c_3^{2m_3-1})+\nu(c_3^{2m_3-1})=2N-1.
\]
In fact, the above argument holds for the other tuple.
Proposition~\ref{prop:splitting} implies
\[
 0<\dim\CC_{2N-1}(E,c_3^{2m_3-1})
   =k_{\nu(c_3^{2m_3-1})}(c_3^{2m_3-1})\leq K_{\nu(c_3^{2m_3-1})}(c_3^{2m_3-1}).
\]
Propositions~\ref{prop:endpoints} and
\ref{prop:characteristic-iteration} now yield
\begin{equation}\label{eq:s4-prime-top}
 K_{\nu(c_3)}(c_3)=1,\qquad
 K_s(c_3)=0\quad(s\ne\nu(c_3)).
\end{equation}

\smallskip
{\bf Step 2.2:} \emph{The normal form of $P_{c_3}$.}

From $G_3=1$ and \eqref{eq:dimension}, we obtain
\[
 p_{3,+}=i(c_3)+p_{3,-}-1,\qquad
 i(c_3)+2p_{3,-}+p_{3,0}\leq4.
\]
Since $i(c_3)\geq3$, there are two possible cases that
\[
 (i(c_3),p_{3,+},p_{3,-})=(4,3,0)
 \quad\text{or}\quad(3,2,0).
\]
In the first case, the normal form is
$P_{c_3}\simeq N_1(1,-1)^{\diamond3}$. Here and below, $\simeq$ denotes a homotopy within
$\Omega^0(P_{c_3})$. Then \(i(c_3^m) = 4m\), \(\nu(c_3^m) = 3\) and \(\hat i(c_3) = 4\) for every \(m\in \mathbb N\). $|2m_3 \hat i(c_3)-2N|=|8m_3-2N|<1$ implies $8m_3=2N$. Thus
\[
 i(c_3^{2m_3})+\nu(c_3^{2m_3})=2N+3.
\]
The nullity is constant along iterations, so
\eqref{eq:s4-prime-top} gives
$\CC_{2N+1}(E,c_3^{2m_3})=0$. All other iterations have already
been excluded from this degree. Hence $M_{2N+1}=0<b_{2N+1}=1$,
which contradicts the existence of the required carrier in degree \(2N+1\).

For $(i(c_3),p_{3,+},p_{3,-})=(3,2,0)$, besides $N_1(1,-1)^{\diamond2}$, the remaining block has
odd index. The blocks $N_1(-1,b)$ and $D(-2)$ give
$\hat i(c_3)=3$. $D(2)$ and $N_1(1,-1)$ have even
index. A block $N_1(1,1)$ contradicts $p_{3,-}=0$.
If it is $I_2$, then $\nu(c_3^m)=4$ for all $m$ and $\hat i(c_3)=4$.
The argument of the first case again excludes this possibility.
The only remaining possibility is therefore
\begin{equation}\label{eq:s4-C-rotation}
 P_{c_3}\simeq N_1(1,-1)^{\diamond2}\diamond R(2\pi\alpha),\qquad
 i(c_3)=3,
 \qquad \hat i(c_3)=2+2\alpha,
\end{equation}
where $1/2<\alpha<1$.
If $\alpha$ is irrational, then $\nu(c_3^m)=2$ for every $m$.
By \eqref{eq:s4-prime-top} and Proposition~\ref{prop:characteristic-iteration},
\[
 \CC_q(E,c_3^m)=0\qquad\text{if }q\ne i(c_3^m)+2.
\]
For this normal form, $S^+_{P_{c_3}}(1)=0$ and $C(P_{c_3})=1$,
so
\[
 i(c_3^{2m_3})=2N-1+2\Delta_3,\qquad
 i(c_3^{2m_3})+\nu(c_3^{2m_3})=2N+1+2\Delta_3.
\]
By \cite[(3.31)]{DuanEll15}, $\Delta_3+\Delta_3'=1$ and we choose the tuple with $\Delta_3=1$.
The central iteration of \(c_3\) can contribute only in degree $2N+3$, so again
$M_{2N+1}=0<b_{2N+1}=1$. Then none of the three central iterations contributes in degree $2N+1$. We conclude that
\begin{equation}\label{eq:s4-alpha-rational}
 \alpha=\frac pq\in\Q\cap(1/2,1),\qquad
 \gcd(p,q)=1,\qquad q\geq3,\qquad3<\hat i(c_3)<4.
\end{equation}
The iteration formulas are
\begin{equation}\label{eq:s4-C-iteration}
\begin{split}
 i(c_3^m)&=2m+2\ceil{m\alpha}-1,\\
 \nu(c_3^m)&=\begin{cases}2,&q\nmid m,\\4,&q\mid m.\end{cases}
\end{split}
\end{equation}

\smallskip
{\bf Step 2.3:} {\it Morse inequalities and the order by length.}

Since $i(c_3)=3$ and $\nu(c_3)=2$, equation
\eqref{eq:s4-prime-top} gives
\[
 \CC_5(E,c_3)=\Q,\qquad \CC_s(E,c_3)=0\quad(s\ne5).
\]
Moreover, $\ceil{2\alpha}=2$ so $i(c_3^2)=7$, and
$i(c_3^m)\geq9$ for $m\geq3$.
The iterations of $c_1$ and $c_2$ contribute only in odd degrees.
Consequently,
\[
 M_j=0\ (j<3),\qquad M_4=M_6=0.
\]
The Morse inequalities give $M_3\geq b_3=1$ and
$M_5\geq b_5=1$. The Morse inequalities in degrees $4$ and $6$
give
\[
 M_4-M_3\geq b_4-b_3=-1,\qquad
 M_6-M_5+M_4-M_3\geq b_6-b_5+b_4-b_3=-2.
\]
Thus
\begin{equation}\label{eq:s4-low-exact}
 M_3=b_3=1,\qquad M_5=b_5=1.
\end{equation}
Every iterate has index at least $6$. Exactly one of
$c_1$ and $c_2$ therefore has index $3$, and we can relabel them so that it is
$c_1$. Since $c_3$ already accounts for the critical module in degree $5$,
\begin{equation}\label{eq:s4-prime-indices}
 i(c_1)=3,\qquad i(c_2)\geq7.
\end{equation}
The chain for \(J=1\) in Corollary~\ref{cor:ordered} has the unique
carriers $c_1,c_3$ in degrees $3,5$, respectively. Hence
\begin{equation}\label{eq:s4-prime-length-order}
 E(c_1)<E(c_3),\qquad L(c_1)<L(c_3).
\end{equation}

\smallskip
{\bf Step 2.4:} {\it The infinite variational visibility and resonance.}

Next we consider the infinite variational visibility of $c_1$ and $c_3$. Without loss of generality,
choose the tuple for which
\begin{equation}\label{eq:s4-pair-two}
 i(c_1^{2m_1})=2N-3,\qquad i(c_2^{2m_2})=2N+3.
\end{equation}
The divisibility condition on $T(c_3)$ gives $q\mid2m_3$. Therefore
\[
 2m_3\hat i(c_3)=4m_3+\frac{4m_3p}{q}\in2\Z,\qquad
 |2m_3\hat i(c_3)-2N|<1,
\]
and hence $2m_3\hat i(c_3)=2N$. Formula~\eqref{eq:s4-C-iteration} yields
\[
 i(c_3^{2m_3})=4m_3+\frac{4m_3p}{q}-1=2N-1,
 \qquad \nu(c_3^{2m_3})=4.
\]
As in Step~2.1, only this central iteration can contribute in degree $2N+1$.
The Morse inequality and Proposition~\ref{prop:splitting} give
\[
 1=b_{2N+1}\leq M_{2N+1}
 =\dim\CC_{2N+1}(E,c_3^{2m_3})
 =k_2(c_3^{2m_3})\leq K_2(c_3^{2m_3}).
\]
Since $0<2<\nu(c_3^{2m_3})=4$, the endpoint properties imply
\begin{equation}\label{eq:s4-central-bottom-zero}
 K_0(c_3^{2m_3})=K_4(c_3^{2m_3})=0.
\end{equation}
In particular,
\[
 \dim\CC_{2N-1}(E,c_3^{2m_3})
 =k_0(c_3^{2m_3})\leq K_0(c_3^{2m_3})=0.
\]

$c_1^{2m_1-1}$ and $c_2^{2m_2-1}$ are nondegenerate, with indices
\[
 i(c_1^{2m_1-1})=2N-3,\qquad
 i(c_2^{2m_2-1})=2N-i(c_2)\leq2N-7.
\]
By Step~2.1, \(c_3^{2m_3-1}\) contributes only in degree $2N-1$.
Together with \eqref{eq:s4-pair-two},
\eqref{eq:s4-central-bottom-zero} and
Lemma~\ref{lem:localization}, we have that
\begin{enumerate}[label=(\arabic*)]
    \item for degree \(2N-3\), the possible carrier is \(c_1^{2m_1-1}\) or \(c_1^{2m_1}\);
    \item for degree \(2N-1\), the possible carrier is \(c_3^{2m_3-1}\);
    \item for degree \(2N+1\), the possible carrier is \(c_3^{2m_3}\).
\end{enumerate}
Applying Theorem~\ref{thm:fr-carriers} in degree
$2N-3=2(N-3)+3$ and $2N+1=2(N-1)+3$, we find that their carriers are iterates of $c_1,c_3$, respectively.
As $N\to\infty$, $c_1$ and $c_3$ are therefore infinitely
variationally visible. Theorem~\ref{thm:visible-resonance} gives
\begin{equation}\label{eq:s4-AC-resonance}
 \hat i(c_1)=\rho L(c_1),\qquad  \hat i(c_3)=\rho L(c_3)
 \qquad\text{for some }\rho>0.
\end{equation}
Combining this with \eqref{eq:s4-prime-length-order}, we obtain
\begin{equation}\label{eq:s4-mean-order}
  \hat i(c_1)< \hat i(c_3).
\end{equation}

\smallskip
{\bf Step 2.5:} {\it A contradiction with the strict ordering.}

Choose the same tuple with
\[
 |\eta_j|<\varepsilon,\qquad
 2\varepsilon<\min\{\hat i(c_3)- \hat i(c_1), \hat i(c_3)\}.
\]
The resonance identity gives
\[
\begin{aligned}
 \rho\bigl(L(c_1^{2m_1-1})-L(c_3^{2m_3-1})\bigr)
 &= \hat i(c_3)- \hat i(c_1)+\eta_1-\eta_3
 >\hat i(c_3)- \hat i(c_1)-2\varepsilon>0,\\
 \rho\bigl(L(c_1^{2m_1})-L(c_3^{2m_3-1})\bigr)
 &=\hat i(c_3)+\eta_1-\eta_3
 >\hat i(c_3)-2\varepsilon>0.
\end{aligned}
\]
Thus the unique carrier in degree $2N-1$ is shorter than both possible
carriers at $2N-3$ in the case (1) above. This contradicts
the strict ordering in \eqref{eq:s4-four-chain}, since
energy and length have the same order.
\end{proof}

\subsection{The non-hyperbolicity result in odd dimensions}

Here we always choose \(\varepsilon>0\) sufficiently small for
Lemma~\ref{lem:localization} and Proposition~\ref{prop:top}.
For a hyperbolic prime closed geodesic \(h\),
Theorem~\ref{thm:long-iteration} and \(\hat i(h)>d\) imply
\[
 i(h^m)=mI,\qquad \nu(h^m)=0,\qquad
 I:=i(h)=\hat i(h)\in\Z,\quad I\geq d+1.
\]
For sufficiently small error, there exist some \(m_h\) such that
\(|2m_hI-2N|<1\) which forcing \(2m_hI=2N\). For every other iteration \(m \neq 2m_h\),
\[
 |i(h^m)-2N|=|m-2m_h|I\geq d+1.
\]
Thus the only possible contribution of a hyperbolic prime closed
geodesic in \([2N-d,2N+d]\) comes from its central iteration and
occurs in degree \(2N\).

\begin{proof}[Proof of Theorem~\ref{thm:odd-stability}]
Let \(d=n-1=2k\) with \(k\geq2\), and let \(p\) denote the number of
non-hyperbolic prime closed geodesics. By
Lemma~\ref{lem:comparison}, every prime closed geodesic \(c\) satisfies
\(i(c)\geq d\) and \(\hat i(c)>d\). Suppose that \(p\leq k\) for a contradiction.

\smallskip
{\bf Step 1:} {\it Counting the non-hyperbolic iterates in degrees greater than \(2N\).}

Choose an arbitrary tuple \((N,m_1,\ldots,m_p)\) satisfying the requirements as above.
First, we consider the \(k-1\) degrees
\begin{equation}\label{eq:odd-positive-fr}
 2N+2,\ 2N+4,\ldots,2N+d-2
\end{equation}
which are strictly between \(2N\) and \(2N+d\). Since hyperbolic geodesics can contribute only in degree \(2N\),
the carriers in these degrees must be central iterates of non-hyperbolic prime closed geodesics.
The strict ordering of the critical values makes
these \(k-1\) prime closed geodesics distinct.

Proposition~\ref{prop:top} provides an additional carrier
\(c_*^{2m_*}\) whose only nonzero critical module is in degree
\(2N+d\). The prime closed geodesic \(c_*\) is non-hyperbolic
and distinct from the preceding \(k-1\) geodesics. Hence \(p\geq k\),
and the assumption gives \(p=k\). The \(k\) geodesics with the central iterations
identified above therefore account for all non-hyperbolic prime
closed geodesics.

\smallskip
{\bf Step 2:} {\it The carrier in degree \(2N\) is hyperbolic.}

By a similar argument, the carrier in degree \(2N\) must be a geodesic with the central iteration,
and the critical module of \(c_*^{2m_*}\) vanishes at \(2N\). Each of
the other \(k-1\) non-hyperbolic geodesics with the central iterations already correspond to
a higher critical value in \eqref{eq:odd-positive-fr}, so it cannot
correspond the strictly smaller value in degree \(2N\).
Thus this carrier in degree \(2N\) is \(h^{2m_h}\) for a hyperbolic prime
closed geodesic \(h\) and some \(m_h>0\).

Take \(N\) sufficiently large such that all indices in this step
and in \eqref{eq:odd-positive-fr} exceed the uniform bound of
Corollary~\ref{cor:visibility-cutoff}. Then the prime closed geodesics
underlying these carriers are infinitely variationally
visible, and their mean indices and lengths have the same positive
ratio \(\rho\). For the hyperbolic carrier, the identity
at the beginning of this section gives
\begin{equation}\label{eq:odd-middle-length}
 2m_h\hat i(h)=2N,\qquad
 L(h^{2m_h})=\frac{2N}{\rho}.
\end{equation}

\smallskip
{\bf Step 3:} {\it Determine the signs of the mean-index errors.}

For each \(c_j^{2m_j}\) in
\eqref{eq:odd-positive-fr}, the strict ordering in this degree
sequence implies
\[
 L(c_j^{2m_j})>L(h^{2m_h}).
\]
Morover, all of them satisfy the resonance identity. By multiplying this inequality by \(\rho\) and using
\eqref{eq:odd-middle-length}, we obtain
\[
 \eta_j=2m_j\hat i(c_j)-2N>0.
\]
For the one remaining non-hyperbolic prime closed geodesic \(c_*\),
Proposition~\ref{prop:top} and Lemma~\ref{lem:band} imply
\[
 2N+d=i(c_*^{2m_*})+\nu(c_*^{2m_*})
 \leq2m_*\hat i(c_*)+d=2N+\eta_*+d.
\]
Thus \(\eta_*\geq0\).
So we have proved that every sufficiently large tuple chosen above must satisfy
\begin{equation}\label{eq:odd-sign-obstruction}
 \begin{aligned}
  \eta_j&\geq0\quad\text{for every non-hyperbolic prime closed geodesic},\\
  \eta_j&>0\quad\text{for at least }k-1\text{ of them}.
 \end{aligned}
\end{equation}
For each such tuple, Step~1 accounts for all \(k\) non-hyperbolic
prime closed geodesics. Hence \eqref{eq:odd-sign-obstruction} applies
to each of them.

If one of the non-hyperbolic prime closed geodesics has irrational
mean index, fix that geodesic. Lemma~\ref{lem:signed-tuple}(ii)
provides arbitrarily large tuples with the same
\(\varepsilon\) and divisibility requirements, and with a negative error
for this fixed geodesic. Steps 1--3 apply to these tuples, which
contradicts the first line of \eqref{eq:odd-sign-obstruction}.
If all the non-hyperbolic mean indices are rational, we can use
Lemma~\ref{lem:signed-tuple}(i) to make all their errors zero
simultaneously. This contradicts the second line of
\eqref{eq:odd-sign-obstruction}, since \(k-1\geq1\).
Both cases are impossible, so \(p\geq k+1\).
\end{proof}

\end{document}